\documentclass[10pt]{article}
 \usepackage{color}
\usepackage{amssymb}
\usepackage{amsthm}
\usepackage{epstopdf}
\usepackage{graphicx}
\usepackage{mathtools}
\usepackage[usenames,dvipsnames, table]{xcolor}
\usepackage[hang]{caption}
\usepackage{hyperref}
\usepackage{lscape} 
\usepackage{oldgerm}  
\usepackage{pgfplotstable}
\usepackage{longtable}
\usepackage{array}   

\newcolumntype{L}{>{$}l<{$}} 

\def\be{\begin{equation}}\def\ee{\end{equation}}
\def\bee{\begin{enumerate}}\def\eee{\end{enumerate}}
\def\bei{\begin{itemize}}\def\eei{\end{itemize}}
\newcommand{\icon}[1]{\includegraphics[height=40pt]{#1}}
\newcommand{\SU}[1]{\mathrm{SU(#1)}}

\newcommand{\nco}{\newcommand}

\nco{\one}{\ensuremath{\,\,\mathrm{l}\!\!\!1}} 
\nco{\NN}{\mathbb{N}}
\nco{\ZZ}{\mathbb{Z}}
\nco{\QQ}{\mathbb{Q}}
\nco{\RR}{\mathbb{R}}
\nco{\CC}{\mathbb{C}}
\nco{\HH}{\mathbb{H}}
\nco{\OO}{\mathbb{O}}

\nco{\red}{\color{red}}
\nco{\blue}{\color{blue}}
\nco{\cyan}{\color{cyan}}
\nco{\brown}{\color{Magenta}}

\nco{\magenta}{\normalcolor}

\nco{\violet}{\color{violet}}
\nco{\redend}{\normalcolor}
\nco{\magentaend}{\normalcolor}

\def\ie{{\it i.e. }}
\def\ommit#1{{}}
\def\({\left(}
\def\){\right)}

\def\ie{{\it i.e.,\/}\ }
\def\ie{{\rm i.e.,\/}\ }

\def\be{\begin{equation}}\def\ee{\end{equation}}
\def\bea{\begin{eqnarray}}\def\eea{\end{eqnarray}}

\nco{\rnc}{\renewcommand}
\rnc{\title}[1]{{\Large\bf\mbox{}\\\medskip#1\bigskip\medskip\\}}
\rnc{\author}[1]{{\large #1\smallskip\\}}
\nco{\address}[1]{{\em #1\medskip\\}}

\def\mydate{November 2019} 
\def\mycontribution{Contribution to the volume \\ Theoretical Physics, Wavelets, Analysis, Genomics\\
An Indisciplinary Tribute to Alex Grossmann\\
(Springer International Publishing)
}

\begin{document}
\begin{titlepage}
\begin{center}
\title{$\SU{3}$ higher roots and their lattices}
\medskip
\author{Robert Coquereaux} 
\address{Aix Marseille Univ, Universit\'e de Toulon, CNRS, CPT, Marseille, France\\
Centre de Physique Th\'eorique}

\bigskip\medskip
\mydate

\vskip 1.cm
\mycontribution
\bigskip\medskip

\begin{abstract}
\noindent {
After recalling the notion of higher roots (or hyper-roots) associated with ``quantum modules''  of type $(G, k)$,  for $G$ a semi-simple Lie group and $k$ a positive integer,  following the definition given by A. Ocneanu in 2000,  we study the theta series of their lattices.
Here we only consider the higher roots associated with quantum modules (aka module-categories over the fusion category defined by the pair $(G,k)$) that are also ``quantum subgroups''.
For $G=\SU{2}$ the notion of higher roots coincides with the usual notion of roots for ADE Dynkin diagrams and the self-fusion restriction (the property of being a quantum subgroup) selects the diagrams of type $A_{r}$, $D_{r}$ with $r$ even, $E_6$ and $E_8$; their theta series are well known.
In this paper we take $G=\SU{3}$,  where the same restriction selects the modules ${\mathcal A}_k$,  ${\mathcal D}_k$  with $mod(k,3)=0$, and the three exceptional cases ${\mathcal E}_5$, ${\mathcal E}_9$ and ${\mathcal  E}_{21}$. 
The theta series for their associated lattices are expressed in terms of modular forms twisted by appropriate Dirichlet characters. 
}
\end{abstract}
\end{center}

\vspace*{70mm}
\end{titlepage}

\section{Introduction}

Root systems of Lie algebras are higher root systems of type $G=\SU{2}$ and  generate  lattices whose properties and associated theta series are well-known.
Here we consider higher root systems of type $G=\SU{3}$ using a general definition given by A. Ocneanu in 2000 \cite{Ocneanu:Bariloche}, see also the more recent reference \cite{Ocneanu:Harvard_roots}.
Such systems are classified by ``quantum modules'' or ``quantum subgroups'' \ie using a categorical language, by module-categories over the modular fusion category defined by a pair $(Lie(G),k)$, where $k$, the level, is a non-negative integer. 
Quantum modules are characterized by graphs, that, for $\SU{3}$, have been obtained long ago in the framework of conformal field theories {\empty  (see sec.~17.10 of the book} \cite{YellowBook},  \cite{DiFrancescoZuber}, \cite{Ocneanu:Bariloche}\footnote{This reference stresses the importance of a cocycle condition for triangular cells, a condition that is explicitly worked out in \cite{CIG:TriangularCells}, see also \cite{EvansPughSU3}.}),
and generalize the usual ADE Dynkin diagrams.
In the present paper we restrict our attention to an even smaller family whose members are called ``quantum subgroups'', and sometimes nicknamed ``exceptional modules with self-fusion'' (they can be defined as connected \'etale algebras in a modular fusion category).
To every quantum module is associated a higher root system. Such systems generate lattices, and our main purpose, after having recalled in the first two sections the necessary definitions and concepts, is to provide closed expressions for the corresponding theta series.
By exhibiting Gram matrices encoding the geometry of lattices we express their theta series in terms of modular forms twisted by appropriate Dirichlet characters. More information about those lattices can be found in \cite{RC:hyperroots}, where ``higher roots'' are called ``hyper-roots''. 
The first few terms of several series had already been obtained by A. Ocneanu \cite{Ocneanu:Bariloche},  \cite{Ocneanu:posters}, and our expressions, that can be expanded to arbitrary orders,  agree with these older results.
Using explicit Gram matrices leads to a technique that is only suitable for low levels but we hope that this work will trigger new developments and  pave the way for more general results.

I learned a lot of group theoretical concepts from Alex, mostly in problems related to quantization. He also introduced me to the fascinating world of elliptic functions and modular forms --- it was in a cosmological framework, {\empty some time before the wavelet period started}.
This contribution to Alex memory has a little to do with quantization, very little to do, as far as I can see, with the theory of wavelets, and probably nothing to do with cosmology (but who knows ?).
Nevertheless, I am sure that Alex would have liked the mixture of geometry and analysis involved in the present work, even if its application to physical sciences still belongs to the realm of science-fiction, and, in any case, Alex was such a kind person that he would never had told me the opposite!

\paragraph{Notations.} {\empty To dissipate possible misunderstandings, the neophyte should maybe read section 2.1 first}.
Warning: the subindex $k$ used for {\sl script} symbols stands for the chosen level, it does not denote the number $r$ of simple objects (the rank). 
Calling $g^\vee$ the dual Coxeter number of $G$, the number $N=g^\vee+k$ is the ``generalized Coxeter number'' or ``altitude'' of the quantum module ${\mathcal E}_k(G)$; remember that $g^\vee=2$ for $G=\SU{2}$ and $g^\vee=3$ for $G=\SU{3}$. 
In the former case,  $N$ coincides with the dual Coxeter number of the Lie algebra defined by the associated Dynkin diagram (quantum McKay correspondence).\\
Quantum subgroups are as follows:

For $G=\SU{2}$, we have ${\mathcal A}_k(\SU{2}) = A_{k+1}$ with $r=k+1$, ${\mathcal D}_{k=4s}(\SU{2}) = D_{2s+2}$ with $r=2s+2$, ${\mathcal E}_{10}(\SU{2}) = E_{6}$ with $r=6$, ${\mathcal E}_{28}(\SU{2}) = E_{8}$ with $r=8$.
The associated graphs (Dynkin diagrams) are well known. {\empty Above and below: the symbol $r$ denotes the number of vertices.}

For $G=\SU{3}$ the quantum subgroups are {\empty denoted} ${\mathcal A}_{k}(\SU{3})$ with $r=(k+1)(k+2)/2!$, ${\mathcal D}_{k=3s}(\SU{3})$ with $r=\tfrac{1}{3} (\tfrac{(k+1)(k+2)}{2} -1)+3$),
 ${\mathcal E}_5(\SU{3})$ with $r=12$,  ${\mathcal E}_9(\SU{3})$ with $r=12$, and ${\mathcal E}_{21}(SU(3))$ with $r=24$.
We shall often drop the reference to $SU(3)$ in the above notations since we are mostly interested in that case.
The graphs associated with ${\mathcal A}_{1}$, ${\mathcal A}_{2}$, ${\mathcal A}_{3}$, ${\mathcal A}_{4}$, ${\mathcal D}_{3}$, ${\mathcal D}_{6}$, ${\mathcal E}_5$, ${\mathcal E}_9$, ${\mathcal E}_{21}$ are sketched below\footnote{{\empty Of course they are \underline{not} Dynkin diagrams}}.
\begin{center}
\icon{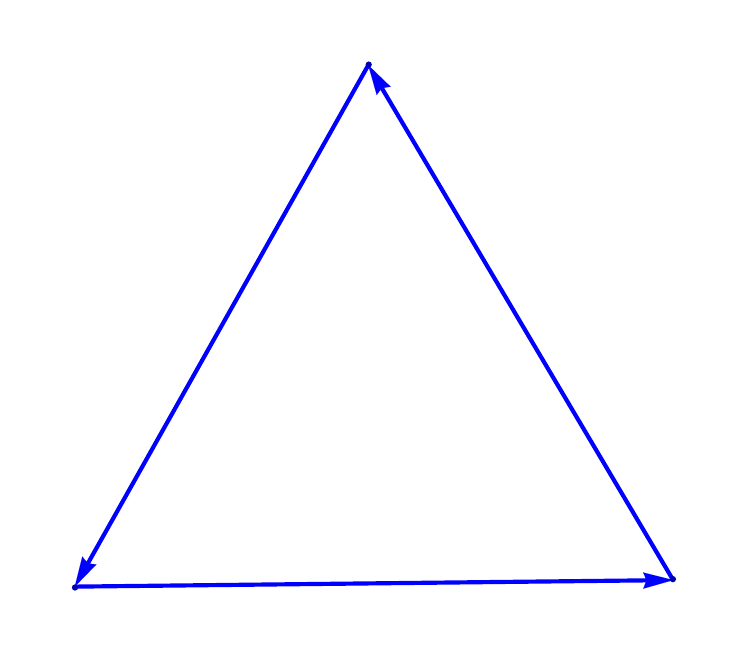} \; \icon{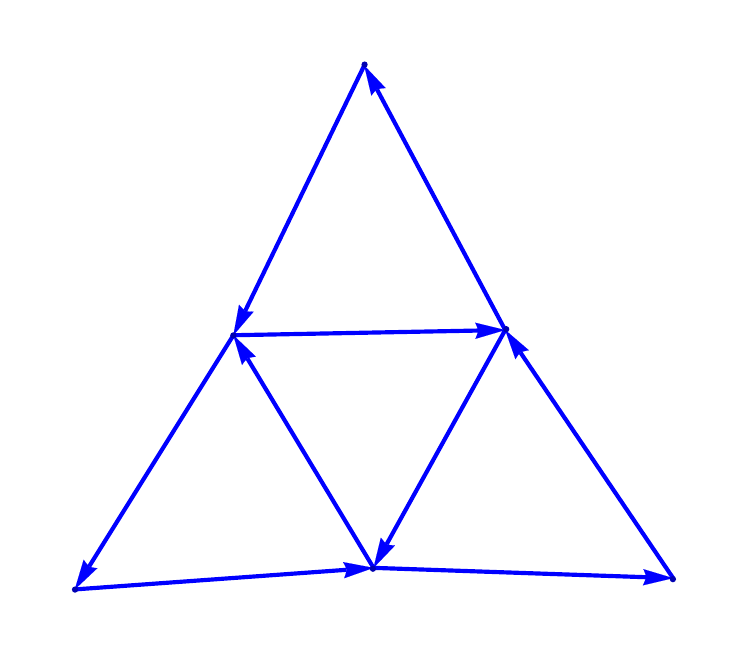} \; \icon{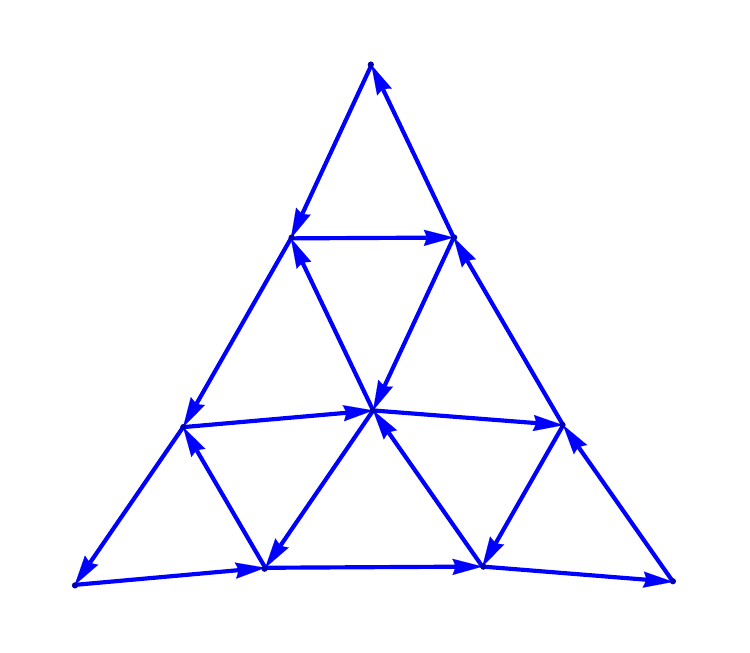} \; \icon{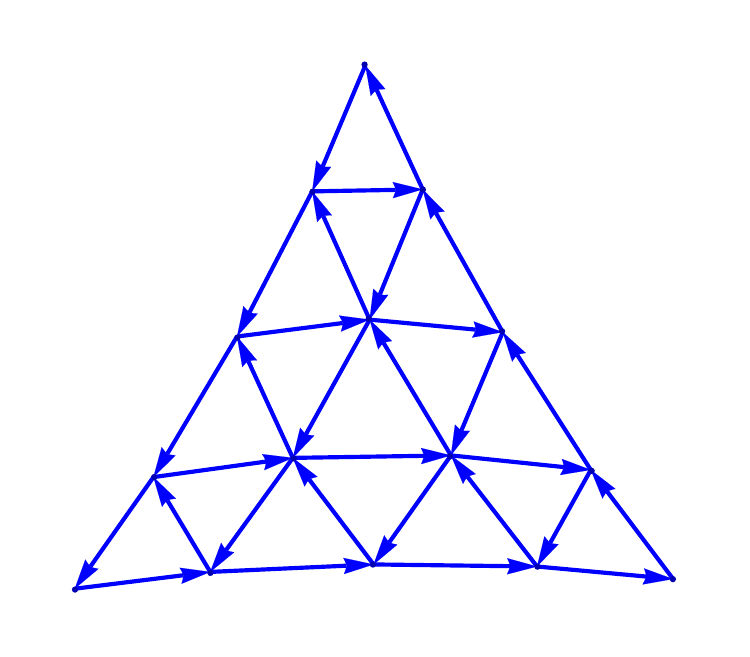}
\icon{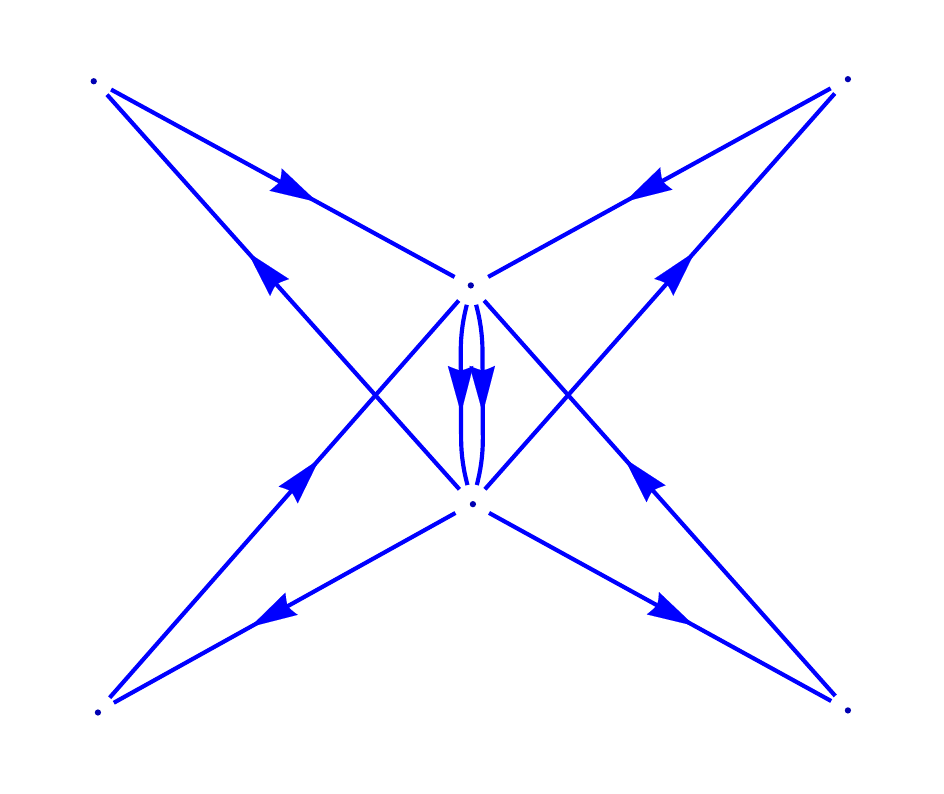} \; \icon{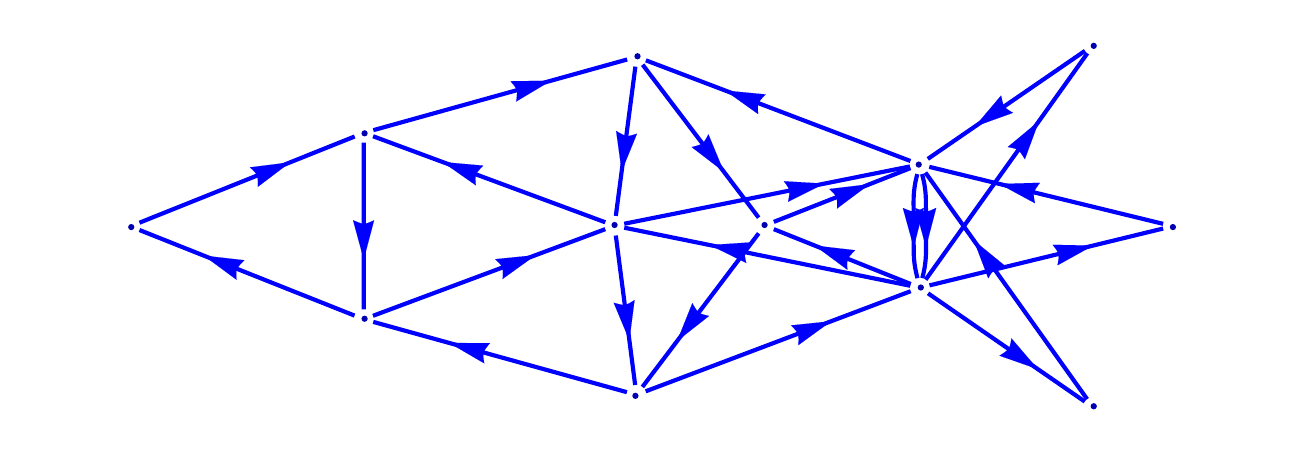}
\icon{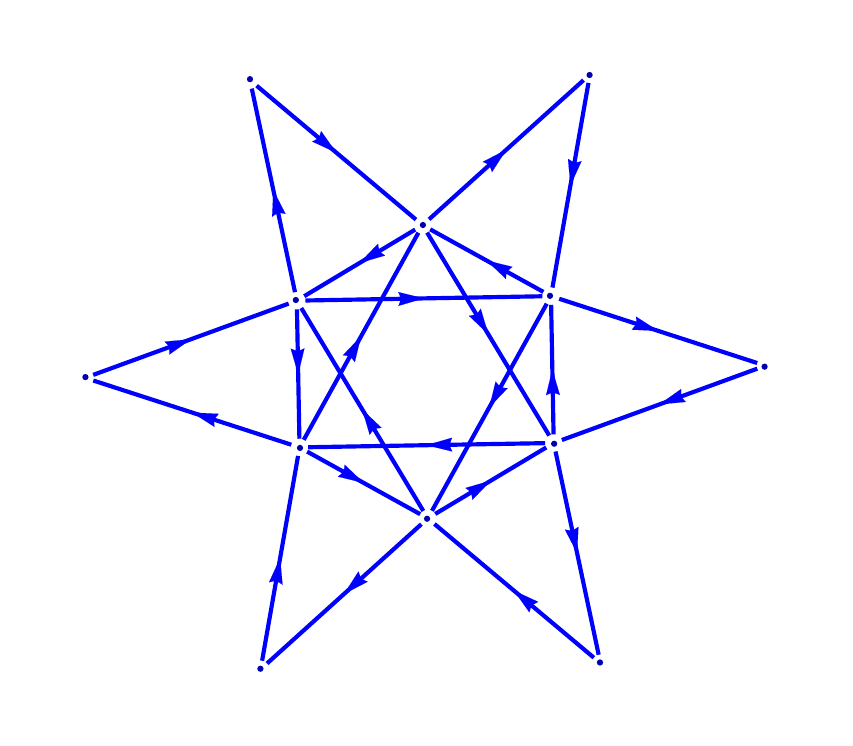} \; \icon{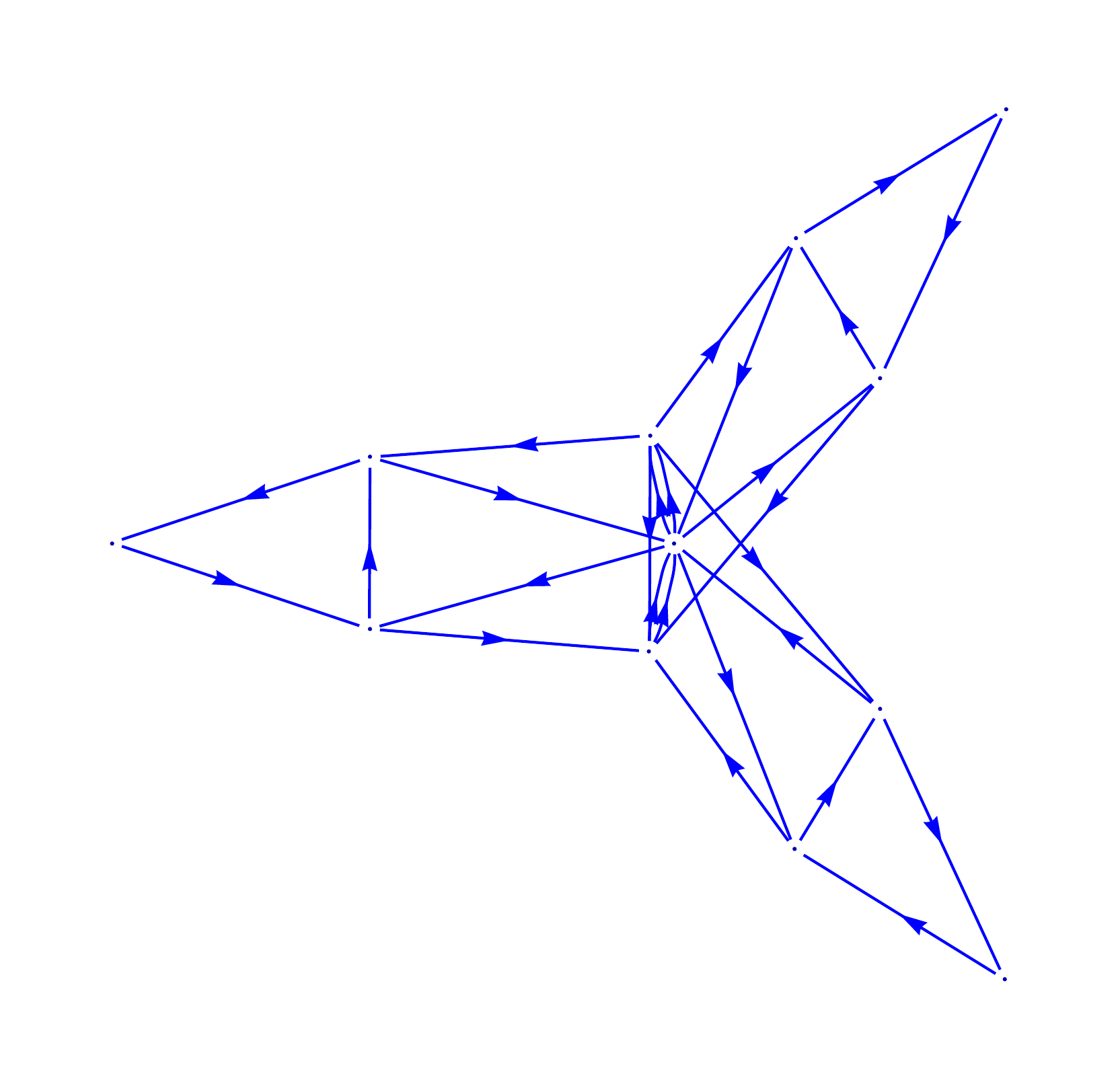} \; \icon{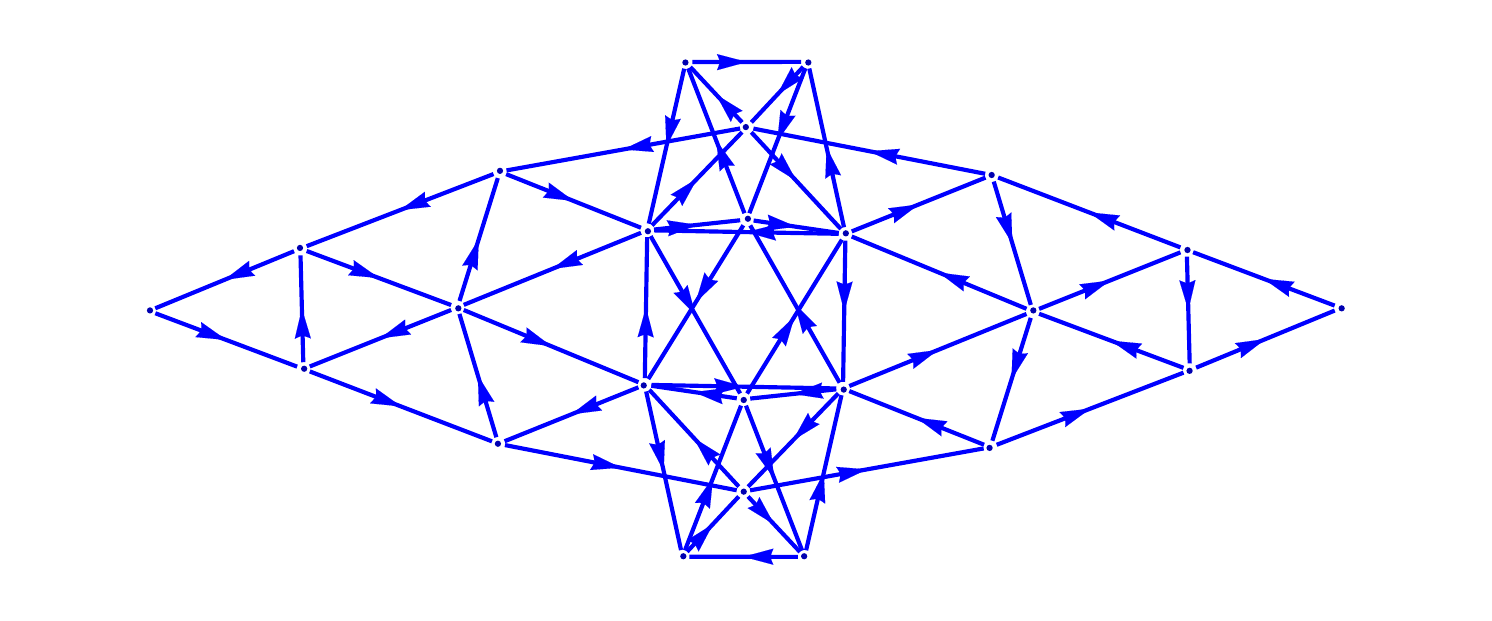}
\end{center}

\paragraph{For the category-minded reader.}
The category ${\mathcal A}_k$ of integrable modules of the affine Kac-Moody algebra associated with $Lie(G)$ at level $k$, see e.g. \cite{Kac:book}, is equivalent (an equivalence of modular tensor categories), 
 \cite{Finkelberg}, \cite{Huang}, \cite{KazhdanLusztig}, to a category constructed in terms of  representations of the quantum group $G_q$ at  root of unity ${\mathfrak q} =  exp(\frac{i \pi}{g^\vee+k})$  ---take the quotient of the category of tilting modules by the additive subcategory generated by indecomposable modules of zero quantum dimension.   A ``quantum module'' denoted ${\mathcal E}_k$  or ${\mathcal E}_k(G)$ is, by definition, a module-category\footnote{It is defined by a monoidal functor from ${\mathcal A}_k$ to the category of endofunctors of an abelian category ${\mathcal E}$, \cite{Ostrik}.  Terminological warning: a module-category is usually not a modular category!} over ${\mathcal A}_{k}={\mathcal A}_{k}(G)$.
 To each module-category over ${\mathcal A}_k$ one associates a system of higher roots. 
  Here we consider only higher roots associated with quantum modules that are also quantum subgroups, \ie module-categories that define a commutative \'etale algebra (also called ${\mathcal E}$)  in ${\mathcal A}_k$  ---every such algebra  giving rise to an indecomposable module-category over ${\mathcal A}_k$ of the particular kind that we have in mind (consider the category of ${\mathcal E}$-modules in ${\mathcal A}_k$). 
 This discussion can be phrased in terms of weak Hopf algebras, see  \cite{Ostrik}.

\section{Lattices of higher roots}

\subsection{Analogies and warnings}
{\empty 
We want to develop an analogy aimed at introducing the subject to the novice reader --- the present section is therefore rather informal and should be skipped by specialists.

Groups may have subgroups, and the space of representations (or of characters) of a subgroup is a module over the ring of representations (or of characters) of the group. The category of representations of a subgroup is a module-category over the category of representations of the corresponding groups: tensor multiplying a representation of a group by a representation of a subgroup gives a representation of the subgroup.

Here we have a similar situation.  The departure point is a category ${\cal A}_k(G)$, where $G$ is some semi-simple Lie group, this category can be defined in many ways, up to equivalence, but we don't even need to know one to proceed. 
It is enough to know that ${\cal A}_k(G)$ is called ``quantum $G$ at level $k$'' 
and that this category is monoidal, meaning that its objects can be multiplied, like the representations of a group. What matters is that we have a ring, which we also called ${\cal A}_k(G)$, that plays the role of a ring of characters (it is indeed the ring of characters for an appropriate structure), and
that this ring may have modules (like the space of representations of a subgroup),  the ring acts on the module and this is called ``fusion''.
Some of those modules have themselves a compatible ring structure,  so  they play the role of  quantum analogs for (space of representations of) subgroups.  Of course, the ring ${\cal A}_k(G)$ is one of its own modules. 

When $G=\SU{2}$, ${\cal A}_k(G)$ is generated by the analog of the spin $1/2$ representation (the fundamental)  but the difference with the classical case is that the number of simpe objects (irreducible representations, or irreps) is finite: there are only $k+1$ of them. The law of composition of spins is the same as usual~: composing a spin $1/2$ with a spin $j$ gives a representation that decomposes as the sum of two irreps, of spins $j-1/2$ and $j+1/2$, unless $j=(k+1)/2$, because in the quantum case at level $k$, the spin $(k+1+1)/2$ does not exist.  The whole multiplication table of spins is determined by the action of the generator ($j=1/2$) and it is encoded by the graph $A_{k+1}$, a truncated line with only $k+1$ vertices : each vertex is an irrep and the result of multiplying the irrep $j$ by $1/2$ are the (two) neighbors of $j$ on the graph. In conformal field theory (WZW models), these vertices are interpreted as ``primary fields'', but we shall not need that. Every module on the latter ring is described by a finite graph that encodes the multiplication (fusion) of basis elements by the generator of ${\cal A}_k(\SU{2})$, and it happens (as discovered first by physicists) that the graphs characterizing those modules are in one-to-one correspondance with the simply-laced Dynkin diagram. Moreover those with ``self-fusion'' (those for which one can define a product compatible with the previous action), i.e. those that are indeed ``quantum subgroups'', have a graph belonging to a smaller family : the $A_k$ themselves, of course, but also the $D_{even}$, and also $E_6$ (for $k=10$), and $E_8$ (for $k=28$).
When $G$ is $\SU{2}$, this is the end of the classification history, which essentially took place ---with another language--- more than $30$ years ago \cite{CIZ}, when the purpose was to classify the modular invariants of WZW models of type $\SU{2}$ (again, we shall not need that).

If one takes $G=\SU{3}$ rather than $\SU{2}$, the graphs that one obtains, that describe again the action of (the quantum analog of) the defining representation of $\SU{3}$ on the different modules of ${\mathcal A}_k(\SU{3})$, have been obtained about $20$ years ago (see \cite{DiFrancescoZuber}, \cite{Ocneanu:Bariloche}). Those with self-fusion are displayed at the end of the introduction. Again, we do not need here to give the precise construction of those modules, or explain how this classification was obtained, but the interested reader could construct them from their graphs, since the adjacency matrix of each of these graphs defines, by definition, the action of the generator. In the $\SU{3}$ case there are actually two generators, the defining representation and its  conjugate, but here graphs are oriented and the action of the conjugated generator is obtained from the given diagrams by reversing the arrows. 

 {\bf Warning.}
Dynkin diagrams are ubiquitous in Mathematics and they are very often introduced when discussing the classification of simple Lie algebras.  It should be stressed that the diagrams obtained for $G=\SU{3}$  are \underline{not} Dynkin diagrams, since the latter all appear when discussing the classification of modules over ${\mathcal A}_k(\SU{2}$)!
As a matter of fact, when presenting a panoramic view of quantum modules for $\SU{2}$ and $\SU{3}$ at level $k$, one does not need to make any reference to (and do not need to use any information from) the general theory of Lie groups. 
For instance, when studying the module described by the graph $E_8$, also denoted ${\mathcal E}_{28}(\SU{2})$ because it is a module over  ${\mathcal A}_{28}(\SU{2}) = A_{29}$, or when studying quantum field theory models associated to this choice, one only needs to use the properties of $\SU{2}$, not those of the exceptional group $E_8$. The latter plays no role in this discussion despite the fact that its Dynkin diagram does appear.
It remains that Dynkin diagrams do have something to do with semi-simple Lie groups! In particular, concepts like roots and weights can be defined from them, and solely from the knowledge of their diagrams.
Since one can ``reverse the machine'' in the $\SU{2}$ case, i.e. start from the classification of modules over ${\mathcal A}_{k}(\SU{2})$ to {\sl discover} the Dynkin diagrams, and from them, define for instance the notion of roots (one can define the roots of $E_8$ without having to define the group $E_8$ itself), it is tempting to consider the possibility of associating roots (or rather ``higher roots'', or ``hyper-roots'') to the graphs describing the classification of modules over ${\mathcal A}_{k}(\SU{3})$, even though there are not Dynkin diagrams.
The interpretation and use of such higher roots in representation theory of usual Lie groups, or in quantum field theory,  is something that is still largely unknown and goes anyway beyond the scope of this article, but it remains that such higher roots have been defined (for any semi-simple $G$, in \cite{Ocneanu:Bariloche}) and they  are still quite mysterious. Like for the usual roots of simple Lie groups, their $Z$-spans define lattices, and the study of these lattices is the study of the present paper.
}

\subsection{On extended fusion matrices and their periods}
\label{fusioncatbasics}

${\mathcal A}_k$ being a monoidal category, with a finite number of simple objects denoted $m,n,p\ldots$, we consider the corresponding Grothendieck ring and its structure coefficients, the so-called fusion coefficients $N_{mnp}$, where  $m \times n = \sum_p N_{mnp} \, p$. They are encoded by fusion matrices $N_m$ with matrix elements $(N_m)_{np} = N_{mnp}$. 
One may consider module-categories ${\mathcal E}$ associated with ${\mathcal A}_k$. 
Of course, one can take for instance ${\mathcal E}= {\mathcal A}_k$. The fusion coefficients $F_{nab}$ characterize the module structure: $n \times a = \sum_b F_{nab} \, b$, where $a,b,\ldots$ denote the simple objects of ${\mathcal E}_k$. They are encoded either by square matrices $F_n$, with matrix elements $(F_n)_{ab}= F_{nab}$, still called fusion matrices, or by the rectangular matrices\footnote{The $F_n$ are sometimes called ``annular matrices'' when ${\mathcal A}_k$ and ${\mathcal E}_k$ are distinct (if they are the same, then $F_n=N_n$), and the $\tau_a$ are sometimes called (by the author) ``essential matrices''.}  $\tau_a = (\tau_a)_{nb}$, with $(\tau_a)_{nb} =  (F_n)_{ab}$. The simple objects of ${\mathcal A}_k$, or irreps, are labelled by the vertices of the Weyl alcove at level $k$.

With $G = \SU{3}$, the simple objects $n$ (irreps) are labelled by pairs $(p,q)$ of non-negative integers with $p+q \leq k$, {\empty so the matrix $F_n=F_{p,q}$ has matrix elements  $(F_{n})_{ab}=(F_{p,q})_{a,b}$}
and the Chebyshev recursion relations of the $G = \SU{2}$ case (composition of spins) are replaced by
\begin{eqnarray}
F_{(p,q)} &=& F_{(1,0)} \, F_{(p-1,q)} - F_{(p-1,q-1)} -F_{(p-2,q+1)} \qquad \qquad \textrm{if} \; q \not= 0
 \nonumber  \label{recursion}\\
F_{(p,0)} &=& F_{(1,0)} \, F_{(p-1,0)} - F_{(p-2,1)} 
\label{recF} \\
F_{(0,q)} &=& (F_{(q,0)})^T 
\nonumber
\end{eqnarray}
$F_{(0,0)}$ is the identity matrix and $F_{(1,0)}$ and $F_{(0,1)}$ are the two generators\footnote{{\empty The adjacency matrices of the graphs given at the end of the introduction are precisely the (unextended) matrices $F_{(1,0)}$}.}. 
Also, $F_{(q,p)}=F_{(p,q)}^T$.

In some applications one sets to zero the fusion matrices whose Dynkin labels do not belong to the chosen Weyl alcove. This is \underline{{\sl not}} what we do here. 
On the contrary, the idea is
to use the same recursion relations to extend the definition of the matrices $F_n$ at level $k$ from the Weyl alcove to
the fundamental Weyl chamber of $G$ (cone of dominant weights) and to use signed reflections with respect to the hyperplanes of the affine Weyl lattice in order to extend their definition to arbitrary arguments $n \in \Lambda$, 
 the weight lattice of $G$.
By so doing one obtains an infinite family of matrices $F_n$ that we still (abusively) call  ``fusion matrices'', and for which we keep the same notations, although their elements can be of both signs.
It is also useful to shift  (translation by the Weyl vector)  the labelling index of the these matrices to the origin of the weight lattice;
in other words, for $n\in \Lambda$, using multi-indices, we set ${\{n\}} = {(n-1)}$, where the use of parenthesis refers to the usual Dynkin labels.
We therefore introduce the ``$\rho$-shifted notation'' ($\rho$ being the Weyl vector), setting $F_{\{p,q\}} = F_{(p-1,q-1)}$, so that  $F_{\{0,0\}}$ is the zero matrix and $F_{\{1,1\}}=F_{(0,0)}$ is the identity (the latter corresponding to the weight with components $(0,0)$ in the Dynkin basis, \ie to the highest weight of the trivial representation).  The following results belong to the folklore:
Setting $N=k+3$ one has  $F_{\{p,q\}} = 0$ whenever $p=0 \; \text{mod} \; N$, $q=0 \; \text{mod} \; N$ or $p+q=0 \; \text{mod} \; N$.
One also gets immediately the following equalities:  $F_{\{p+N,q\}}= (P.F)_{\{p,q\}}$,  $F_{\{p,q+N\}}= (P^2.F)_{\{p,q\}}$ where $P=F_{\{N-2,1\}}$ is a generator of $\ZZ_3$ (with $P^3=1$) acting by rotation on the fusion graph of ${\mathcal A}_k$ and $F_{\{p+3N,q\}}= F_{\{p,q+3N\}} = F_{\{p+N,q+N\}}=F_{\{p,q\}}$.
The sequence $F_{\{p,q\}}$  is periodic of period  $3 N$ in each of the variables $p$ and $q$ but it is completely characterized by the values that it takes in a rhombus  with $N\times N$ vertices; for this reason, this rhombus $D$ will be called periodicity cell, or periodicity rhombus. We have reflection symmetries (with sign) with respect to the lines $\{p\}=0 \; \text{mod} \; N$, $\{q\}=0 \; \text{mod} \; N$ and $\{p+q\}=0 \; \text{mod} \; N$.
The $F$ matrices  labelled by vertices belonging to the Weyl alcove specified by the choice of a non-negative integer $k$
have non-negative integer matrix elements (the alcove is strictly included in the first half of a periodicity rhombus); those with indices belonging to the other half of the inside of the rhombus have non-positive entries, those with vertices belonging to the walls of the Weyl chamber or to the second diagonal  of the rhombus vanish, and the whole structure is periodic. 
The Weyl group action\footnote{This is the shifted Weyl action: $w \cdot n= w (n+\rho) - \rho$ where $\rho$ is the Weyl vector.} on the $\SU3$ lattice is well known.

Since we have extended the definition of the fusion matrices $F_n$ to allow arguments $n$ belonging to the weight lattice we can do the same for the essential matrices $\tau$'s, keeping the same notation: the indices $a,b$ of $(\tau_a)_{nb}$ still refer to simple objects of ${\mathcal E}$ but the index $n$ labels weights of $G$;  the infinite matrices $\tau_a$ can be thought of as rectangular, with columns indexed by the elements of ${\mathcal E}$ (a finite number) and lines indexed by the weights of $\Lambda$, the weight lattice of $G$.

\subsection{The quiver of higher roots (``the Ocneanu ribbon'')}

\subsubsection{The ribbon ${\mathcal R}$}
\label{sec:ribbon}
Given a module-category ${\mathcal E}$ over ${\mathcal A}_k$ (the same notation ${\mathcal E}$ will also denote the set of isomorphisms classes of its simple objects), we defined, 
for every choice of $a \in {\mathcal E}$, an infinite matrix $\tau_a$ which is a periodic, integer-valued, particular function, on $\Lambda \times {\mathcal E} $.
In many cases the definition domain of $\tau_a$ can be further restricted: indeed, there are many modules ${\mathcal E}$ that have a non trivial grading with respect to the center ${\mathcal Z}$ (here  $\ZZ_3$) of the underlying Lie group
(here $G=\SU{3}$);  in those cases, not only the weights of $G$, its irreducible representations, the simple objects of ${\mathcal A}_k$, but also the simple objects of  ${\mathcal E}$, have a well defined grading (denoted $\partial$) with respect to  ${\mathcal Z}$, and the module structure is compatible with this grading: matrix elements of $\tau_a$ in position $(n,b)$ will  automatically vanish if $\partial n + \partial a \neq \partial b$.  Existence of such a non trivial grading occurs for all the module-categories considered in this paper (cases with self-fusion).
The function $\tau_a$ on $\Lambda \times_{\mathcal Z} {\mathcal E}$, is specified 
by the values that it takes on the finite set ${\mathcal R}^\vee = D \times_{\mathcal Z}  {\mathcal E}$ where $D$ is the period parallelotope.
 The set ${\mathcal R}^\vee$,  a finite rectangular table made periodic, may be thought as a closed ribbon\footnote{The terminology ``ribbon'' comes from A. Ocneanu.}. 
 For the cases that we consider,
  the group ${\mathcal Z}$ acts non trivially and ${\mathcal R}^\vee$ has $r_{\mathcal E}  \; \vert{D}\vert / \vert {\mathcal Z} \vert $ elements, where the rank $r_{\mathcal E}$ is the number of simple objects of ${\mathcal E}$,
  and $\vert{D}\vert =  N^{r_G}$.
The elements of ${\mathcal R}^\vee$ will be called restricted (for reasons given below) higher roots of type $G$ defined by the module ${\mathcal E}_k(G)$.
With $G=\SU{3}$ one has $\vert{{\mathcal R}^\vee}\vert  = r_{\mathcal E}  (k+3)^2 /3$.

The choice of a fundamental irrep $\pi$ of $G$, with the constraint that it should exist at level $k$---so that $\pi$ defines a particular non-trivial simple object of ${\mathcal A}_k(G)$---allows one to associate with a graph ${\mathcal E}^\pi$ (but we shall usually drop the reference to $\pi$ in the notation):  it is the graph of multiplication by $\pi$, sometimes called fundamental fusion graph, fundamental representation graph, nimrep graph, or McKay graph associated with $\pi$. If $\pi$ is complex, edges of ${\mathcal E}$ are oriented; it is actually a quiver since it is a directed graph where loops and multiple arrows between two vertices are allowed.
For $\SU{3}$ there are two fundamental irreps, they are complex conjugate to one another and both appear already at level $1$; the associated directed graphs just differ by a global change of orientation (reverse arrows), for definiteness we choose the fundamental irrep with highest weight $(1,0)$ in the Dynkin basis, \ie $\{2,1\}$ in back shifted (\ie $\rho$ shifted) coordinates.

The weight lattice $\Lambda$ can be considered as a directed graph in the usual way, the direct successors of $\{p,q\}$ being the vertices with components $\{p+1,q\}$, $\{p,q-1\}$ and $\{p-1,q+1\}$; 
its restriction to the dominant Weyl chamber can be interpreted as the representation graph of the $(1,0)$-fundamental representation of $\SU{3}$.
 Being essentially a cartesian product of two directed multigraphs, the set ${\mathcal R}^\vee$ becomes a quiver:  the quiver of higher roots, or ``ribbon''.
 
If $\SU{3}$ is replaced by $\SU{2}$, the graph ${\mathcal E}$ has an adjacency matrix which is  twice the identity minus the Cartan matrix of a simple Lie group and the above construction leads to an associated quiver of roots; several examples of this construction (for instance the quiver of $E_6$ with its $72$ vertices) can be found in \cite{RC:periodicquivers}.

For usual roots, \ie $\SU{2}$  higher roots, the opposite of a root is a root. However, for $\SU{3}$ higher roots, one can see, using the definition of the period parallelotope $D$, that if $\alpha\in {\mathcal R}^\vee$, $-\alpha$ does not correspond to any vertex of ${\mathcal R}^\vee$. This feature is not convenient. For all purposes it is useful to  generalize the previous definitions, keeping ${\mathcal R}={\mathcal R}^\vee$ for $\SU{2}$ but  setting ${\mathcal R}= {\mathcal R}^\vee  \cup - {\mathcal R}^\vee$
for $\SU{3}$, then  $\vert {\mathcal R} \vert =2  \vert  {\mathcal R}^\vee  \vert$.
The opposite of a higher root (an element of ${\mathcal R}$) is then always a higher root. For higher roots associated with quantum subgroups ${\mathcal E}$ of $\SU{3}$ we have therefore
\be  \vert {\mathcal R}\vert = 2 \vert {\mathcal R}^\vee \vert = \frac{2}{3}  r_{\mathcal E} \, (k+3)^2   \ee
This can be seen as a generalization of the Kostant relation: for $\SU{2}$ higher roots, \ie usual roots, one has  $\vert {\mathcal R}\vert =  r_{\mathcal E} \, g$ where $g=k+2$ is the Coxeter number, and $r_{\mathcal E}$ the number of nodes, of the Dynkin diagram  ${\mathcal E}$.

For $\SU{3}$ higher roots, if one chooses ${\mathcal E} = {\mathcal A}_k$ one has $r_{\mathcal E} = (N-2)(N-1)/2$, with $N=k+3$, then\footnote{The number of higher roots for ${\mathcal A}_k (\SU{3})$ is therefore also given by the number of 
3-cycles in the rook graph $M_N=K_N  \square K_N$, the Cartesian square of the complete graph  $K_N$ on $N$ vertices (one recognizes the A288961 sequence of the OEIS~\cite{OEIS}).}.
\be  \vert {\mathcal R}\vert =  2\vert {\mathcal R}^\vee \vert = (N-2)(N-1) N^2 /3   \ee
From now on, we shall usually not mention ${\mathcal R}^\vee$, the set of restricted higher roots, since ${\mathcal R}$ will be used most of the time.

\subsection{An Euclidean structure on the space of higher roots}
\label{EuclidanStructureHigherRoots}

The inner product of two higher roots $\alpha = (m,a)$ and $\beta=(n,b)$ is defined as
\be
<\alpha, \beta > =  \sum_{w \in {\mathcal W}} \epsilon(w) F_{m-n + w\rho - \rho, a,b} 
\label{scalarproductofhyperroots1}
\ee
where $W$ is the Weyl group.
The above expression generalizes the one obtained\footnote{This was recognized and generalized in \cite{Ocneanu:Bariloche} but it is already present in \cite{Dorey:CoxeterElement}.} for $G=\SU{2}$. 
Another possibility (the approach followed in \cite{Ocneanu:Bariloche} ) is to start from the notion of harmonicity:
The ${\mathcal R}^\vee$ quiver being essentially a cartesian product of two multigraphs, we have a natural notion of harmonicity for the functions defined on its underlying set; 
a point $\alpha$ of ${\mathcal R}^\vee$ specifies a harmonic function, also denoted $\alpha$, defined as the orthonormal projection of the Dirac measure $\delta_\alpha \in {\CC^{\vert {\mathcal R}^\vee \vert}}$ on the subspace
of harmonic functions; one proves that its value on $\beta \in {\mathcal R}^\vee$, is precisely given by Eq.~\ref{scalarproductofhyperroots1}. More generally, higher weights are $\ZZ$-valued functions that are harmonic on the ribbon.
For us it will be enough to take Eq.~\ref{scalarproductofhyperroots1} as a definition of the inner product between higher roots and extend $\langle  \; , \;  \rangle$ by linearity to their linear span.  
One checks that it defines a positive definite\footnote{Using Eq.~\ref{scalarproductofhyperroots1}
one could define a periodic inner product  on  $\Lambda \times_{\mathcal Z}{\mathcal E}$ that would not be positive definite because of the periodicity, 
but we consider directly its non-degenerate quotient, naturally defined on the ribbon $D \times_{\mathcal Z}{\mathcal E}$.}
 inner product and therefore an Euclidean structure on the space of higher roots. This Euclidean space will be denoted~${\mathfrak C}$.
We did not introduce here any higher analogue of the non-simply laced condition, so that higher roots have only one possible length: $<\alpha, \alpha > =  \vert{\mathcal W}\vert$, for all higher roots $\alpha$.
With $G=\SU{3}$, the Weyl group is ${\mathcal S}_3$,  all higher roots have norm $6$ and  the inner product of two higher roots is obtained from Eq.~\ref{scalarproductofhyperroots1}
as the sum of six fusion coefficients\footnote{With $G=\SU{2}$, one recovers the fact that roots of simply laced Dynkin diagrams have norm $2$ and that the inner product of two roots can be written as a sum of two fusion coefficients.}.
Writing $\alpha=(m,a) = ((m_1,m_2),a)$, $\beta=(n,b) = ((n_1,n_2),b)$, setting $\lambda_1=m_1-n_1$, $\lambda_2=m_2-n_2$, and using shifted labels, we obtain
\begin{eqnarray}
& & <\alpha, \beta > \quad =\\
& & \left(F_{\{\lambda_1+1,\lambda_2+1\}}  +  F_{\{\lambda_1-2,\lambda_2+1\}}+F_{\{\lambda_1+1,\lambda_2-2)\}}-F_{\{\lambda_1-1,\lambda_2-1\}}-F_{\{\lambda_1-1,\lambda_2+2\}}-F_{\{\lambda_1+2,\lambda_2-1\}} \right)_{(a,b)}
\label{scalarproductofsu3roots}
\end{eqnarray}

Elements of the higher root lattice, the $\ZZ$--span of higher roots, are called ``higher-root vectors'', and the elements of its dual lattice are ``higher-weight vectors''.

\subsection{Rank of the system}

The dimension $\mathfrak{r} = \text{dim} \; {\mathfrak C}$ of the space of higher roots associated with ${\mathcal E}_k(G)$, in those cases
 where the center ${\mathcal Z}$ acts non-trivially on the set of vertices of ${\mathcal E}_k(G)$, is\footnote{This general result was claimed in the last two slides of \cite{Ocneanu:MSRI} and it can be explicitly checked in all the cases that we consider below.}
$\mathfrak{r} = \text{dim} \; {\mathfrak C} = \frac{ r_{\mathcal E} \, \vert {\mathcal W} \vert }{\vert {\mathcal Z} \vert}$
 where ${\mathcal W}$ is the Weyl group associated with the  simple Lie group $G$.
 
The term ${\vert {\mathcal W} \vert }/{\vert {{\mathcal Z}} \vert }$ cancels out for $G=\SU{2}$ since $ {\mathcal W}$ and $ {\mathcal Z}$ are both isomorphic with $\ZZ_2$ and one then recovers the rank $\mathfrak{r}  = r_{\mathcal E}$ given by the number of vertices of the chosen Dynkin diagram.

For $G=\SU{3}$, ${\mathcal W} ={\mathcal S}_3$, $\vert {\mathcal W} \vert=3!$, and for modules ${\mathcal E}$ with non trivial triality  we have
${\mathcal Z}=Z_3$, therefore $\mathfrak{r}=  2\, r_{\mathcal E}$.
Moreover, if one chooses ${\mathcal E} = {\mathcal A}_k$,  then $\mathfrak{r} = (N-2)(N-1)$ with $N=k+3$.

 Warning: 
 For higher roots of type $\SU{3}$, in contrast with the case of usual roots, one should remember that  $\mathfrak{r}$ and $r_{\mathcal E}$ are related by a factor $2$,  and this implies that the combinatorial structure encoding the inner product in the space of higher roots will require {\sl two copies\/} of the fusion graph ${\mathcal E}$ generalizing the Dynkin diagram.
A naive generalization of the equation $A = 2 - F_{(1)}$  relating the adjacency matrix of Dynkin diagrams to the Cartan matrix could suggest, in the case of $\SU{3}$ higher roots,  to replace $A$ by $6 - (F_{(1,0)}+F_{(0,1)})$, 
some properties of this last matrix and of its inverse are actually investigated in one section of \cite{CoquereauxZuberNuclPhys}, see also \cite{CoquereauxSchieberJMP}, but the lattices obtained from this naive choice  {\sl are not} the lattices of higher roots considered in the present paper.

\subsection{A graphical summary}

\begin{figure}[ht]
\centering
\includegraphics[width=26pc]{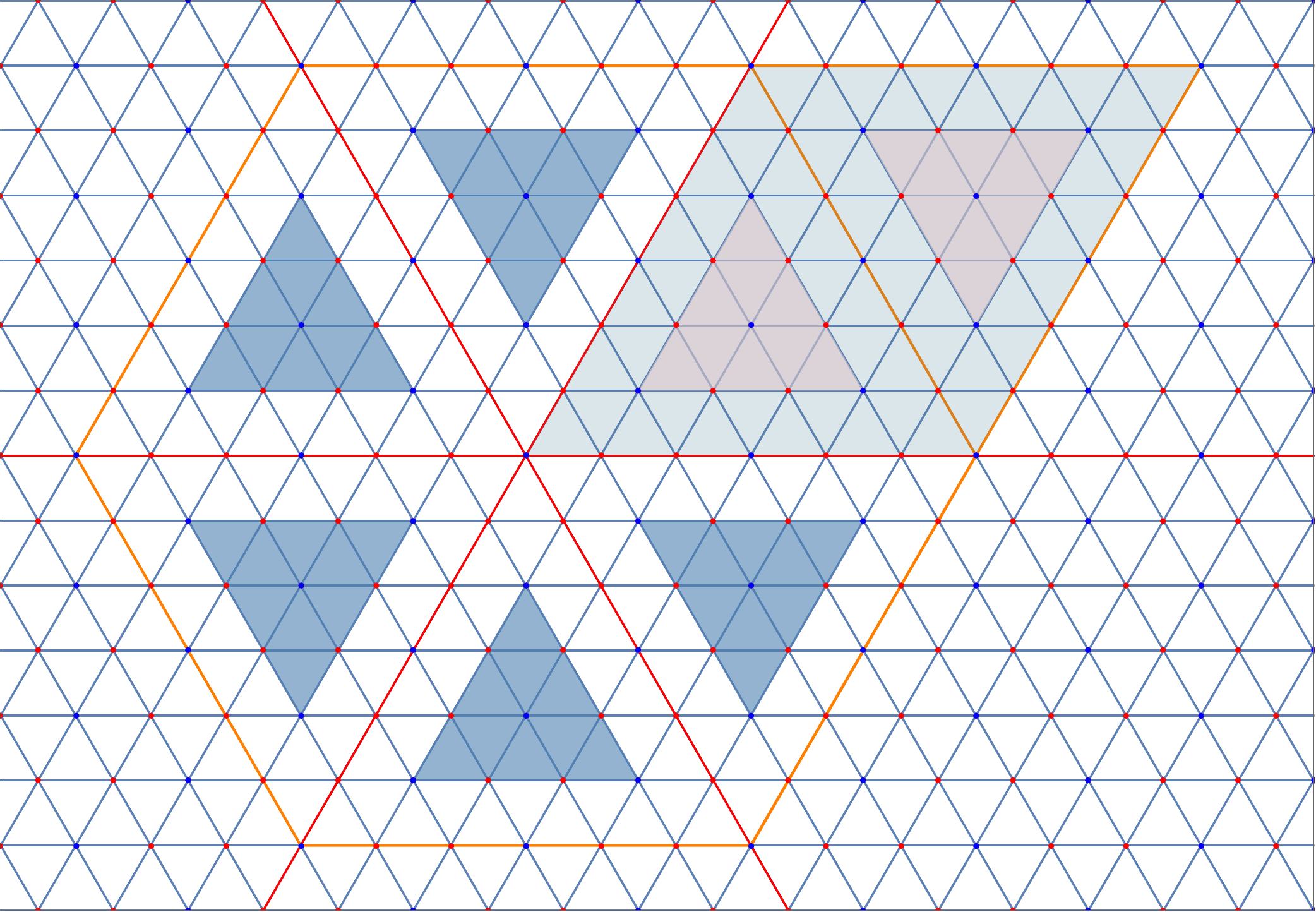}
\caption{\label{hexagonalBackground} The $\SU{3}$ weight lattice $\Lambda$ and the period parallelotope $D$ at level $3$.}
\end{figure}

For illustration  we consider the module ${\mathcal E} = {\mathcal A}_3={\mathcal A}_3(\SU{3})$. 
The weight lattice $\Lambda$ of $\SU{3}$ and the period parallelotope at level $k=3$, a rhombus that we called $D$, are displayed in Fig.~\ref{hexagonalBackground}. 
Weights  are blue or red dots, roots are blue dots.
The edges of $D$ have length $N=k+3=6$, the generalized Coxeter number (altitude). 

\begin{figure}[ht]
\centering
\includegraphics[width=26pc]{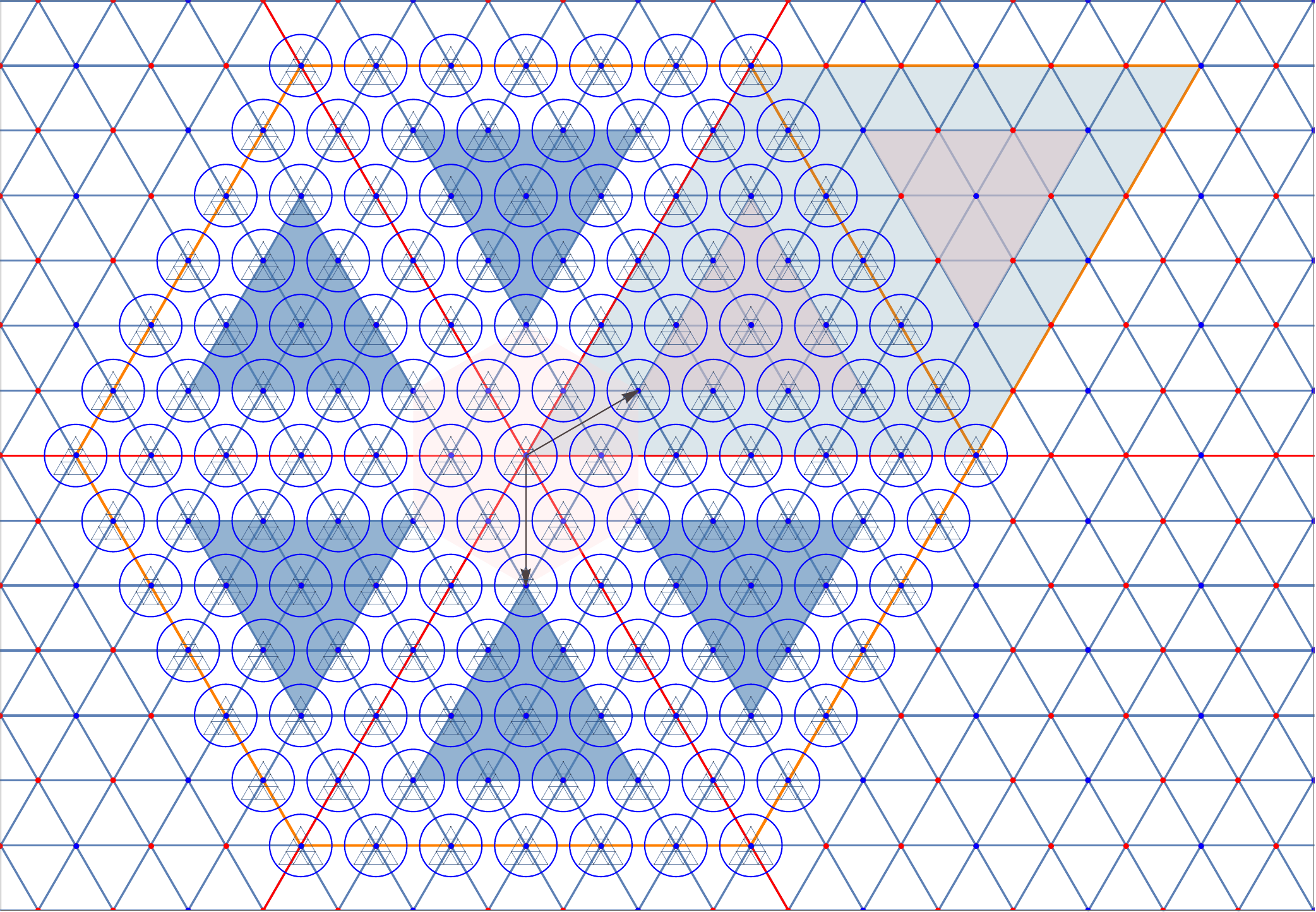}
\caption{\label{ProductOfhexagonalBackgroundByA3} Hexagonal window in the Cartesian product $\Lambda \times {\mathcal A}_3$}
\end{figure}

The next step is to consider the cartesian product $\Lambda \times {\mathcal E}$, which is displayed in Fig.~\ref{ProductOfhexagonalBackgroundByA3}: here ${\mathcal E}$ denotes the graph encoding the module structure (its adjacency matrix, for the chosen example, is the fundamental fusion matrix of $\SU{3}$ at level $3$). So, to each weight one associates a copy of the graph ${\mathcal A}_3$  (we only draw the graphs associated with weights belonging to a rhombus of side $N$).

Each of the  ``small graphs ${\mathcal A}_3$''  have $10$ vertices but not all vertices of those small graphs will represent higher roots since one has to quotient the previous Cartesian product by the action of the center ${\mathcal Z}=Z_3$:   only the large dots of Fig.~\ref{ExtendedParallelogramA3} (Right) define higher roots (the small graphs contain either $3$ of $4$ such large dots, leading to a total of $10\times6^2 / 3 = 120$ higher roots). Because of periodicity one should take only once the small graphs associated with weights belonging to the edges of the period parallelogram, in other words, the top and right edges, in Fig.~\ref{ExtendedParallelogramA3}, do not contribute to the counting, so  that the total is again $(4 \times 2 + 3 \times 4)\times 6 = 120$,  as expected.
{\sl The reader should use a screen display with a large enough magnification factor since, obviously, the described features cannot be seen in a printed version}.
In order not to clutter Figs~\ref{ExtendedParallelogramA3}  and \ref{pathLabelsAlpha23_larger},  the arrows of the resulting periodic quiver are not drawn ---the interested reader may however look at Fig.~2 of reference \cite{RC:hyperroots} where a simpler example is completely described.

\begin{figure}[ht]
\centering
\includegraphics[width=15pc]{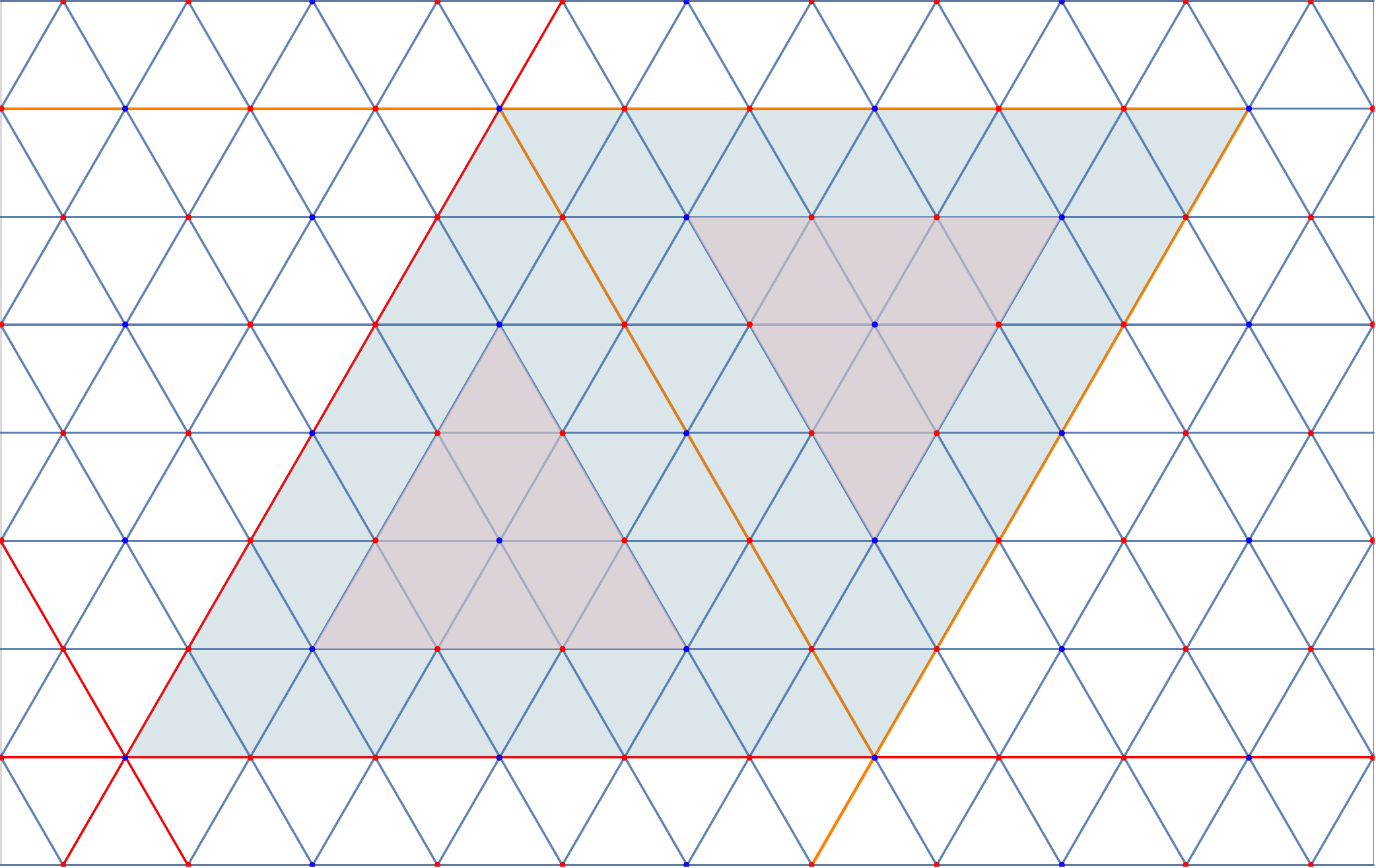}
\quad
\includegraphics[width=15pc, height=3.7cm,trim= 0 -3.3cm 3.3cm -0cm]{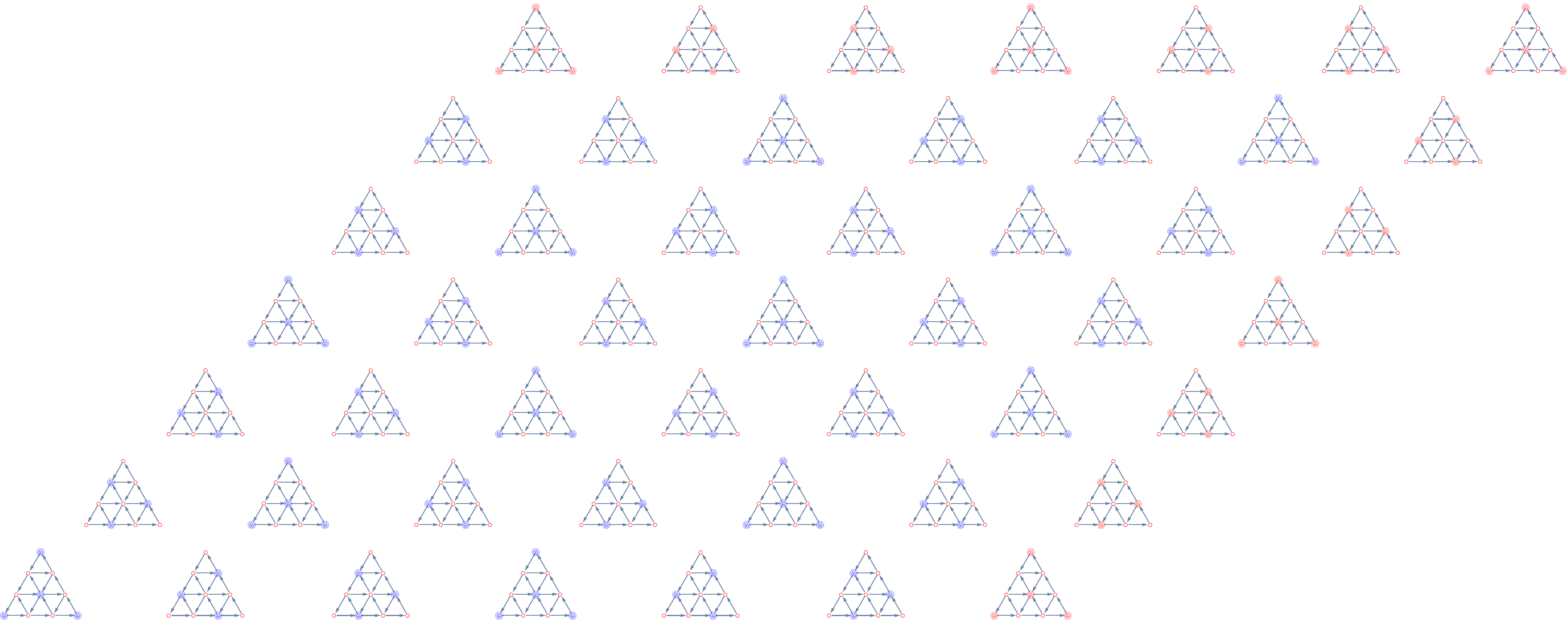}
\caption{\label{ExtendedParallelogramA3} Left: the period parallelogram $D$ of ${\mathcal A}_3(\SU{3})$; its horizontal base has length $N=3+3=6$.
Right: the positions of the $120$ higher roots (blue dots) obtained by taking the Cartesian product of $D$ by the ``small'' fusion graph of ${\mathcal A}_3$ over the center $Z_3$.
Note: $(N+1)^2=7^2$ small graphs are displayed on top of $D$ but it is enough to consider $6^2$ small graphs (with blue dots) since  those located on the top and right edges of $D$ (with red dots) can be obtained using periodicity.}
\end{figure}

Fig.~\ref{pathLabelsAlpha23_larger}  displays the family of inner products between some chosen higher root (the one marked  "$6$" in the picture) and all the higher roots  of the quiver ${\mathcal R}^\vee$ of ${\mathcal A}_3(\SU{3})$.  Those inner products are calculated using Eq.~\ref{scalarproductofsu3roots}.  
\begin{figure}[ht]
\centering
\includegraphics[width=28pc]{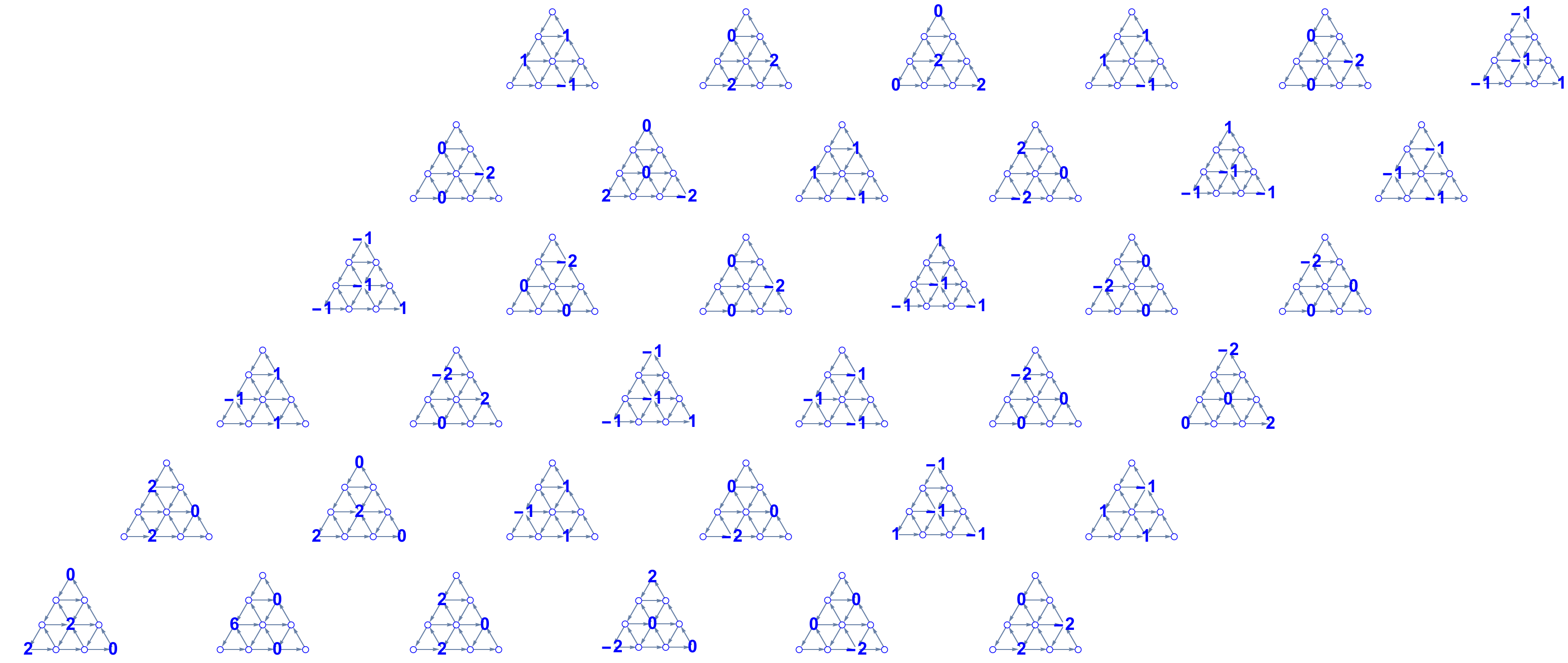}
\caption{\label{pathLabelsAlpha23_larger} The family of inner products between the higher root marked $6$ and all the higher roots of the quiver ${\mathcal R}^\vee$ of ${\mathcal A}_3$ (directed edges are not displayed).}
\end{figure}

\begin{figure}[ht]
\centering
\includegraphics[width=26pc]{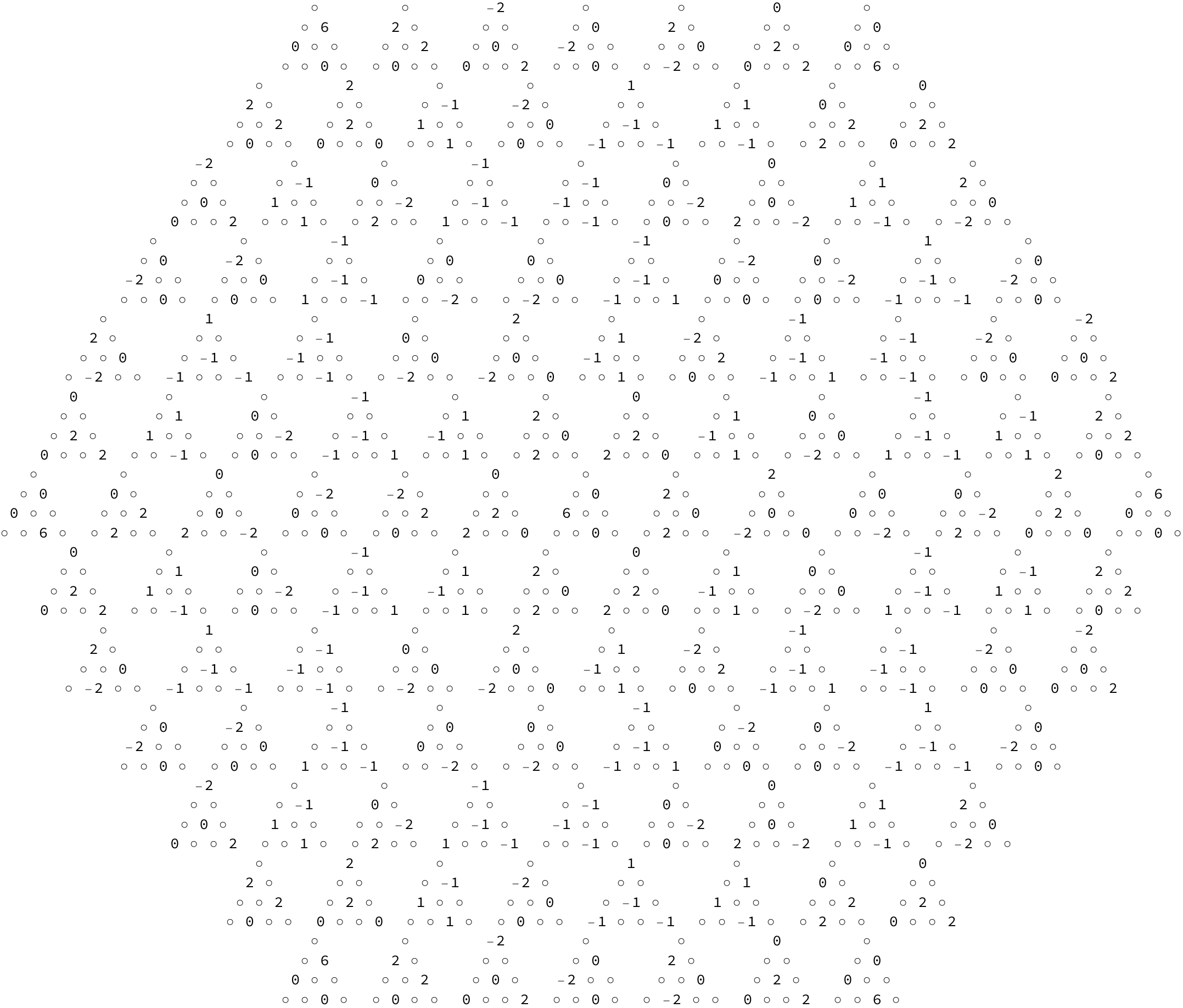}
\caption{\label{CenteredHexagonAlpha23} ${\mathcal A}_3(\SU{3})$: hexagon displaying the periodicity of the scalar product between a chosen higher root (located at the center of the hexagon) and all the higher roots. Note: the chosen higher root is the same as the one chosen in Figs.~\ref{pathLabelsAlpha23_larger} or \ref{A3scalarproductssimplerootsA3}, but the whole picture has been shifted to the right by one unit, in order for this higher root to be at the center of the hexagon.}
\end{figure}

One obtains in this way a periodic function which is completely characterized by the values that it takes in ${\mathcal R}^\vee$, but, as a function on a discrete subset of $\Lambda \times {\mathcal E}$, it is periodic of period $3N$ for each of the directions parallel to the walls of the Weyl lattice (see our discussion in Sec.~\ref{fusioncatbasics}); one can also describe the resulting periodic function by displaying the values that  it takes on the hexagon given in Fig.~\ref{CenteredHexagonAlpha23} ---here we do not  even draw the edges of the small graphs.
The same hexagon can be useful to explicitly check the harmonic property mentioned in Sec.~\ref{EuclidanStructureHigherRoots}: the sum of the values taken by a given higher root on its preimage relative to the ``large directed graph'' $\Lambda$ is equal to the sum of the values that it takes on its preimage relative to the ``small directed graph'' ${\mathcal E}$ (see Fig.~4 of  \cite{RC:hyperroots}  for illustration).

\ommit{
\begin{figure}[ht]
\centering
\includegraphics[width=12pc]{A3parallelogramOfRootsEmpty.pdf}
\includegraphics[width=12pc, height=2.5cm]{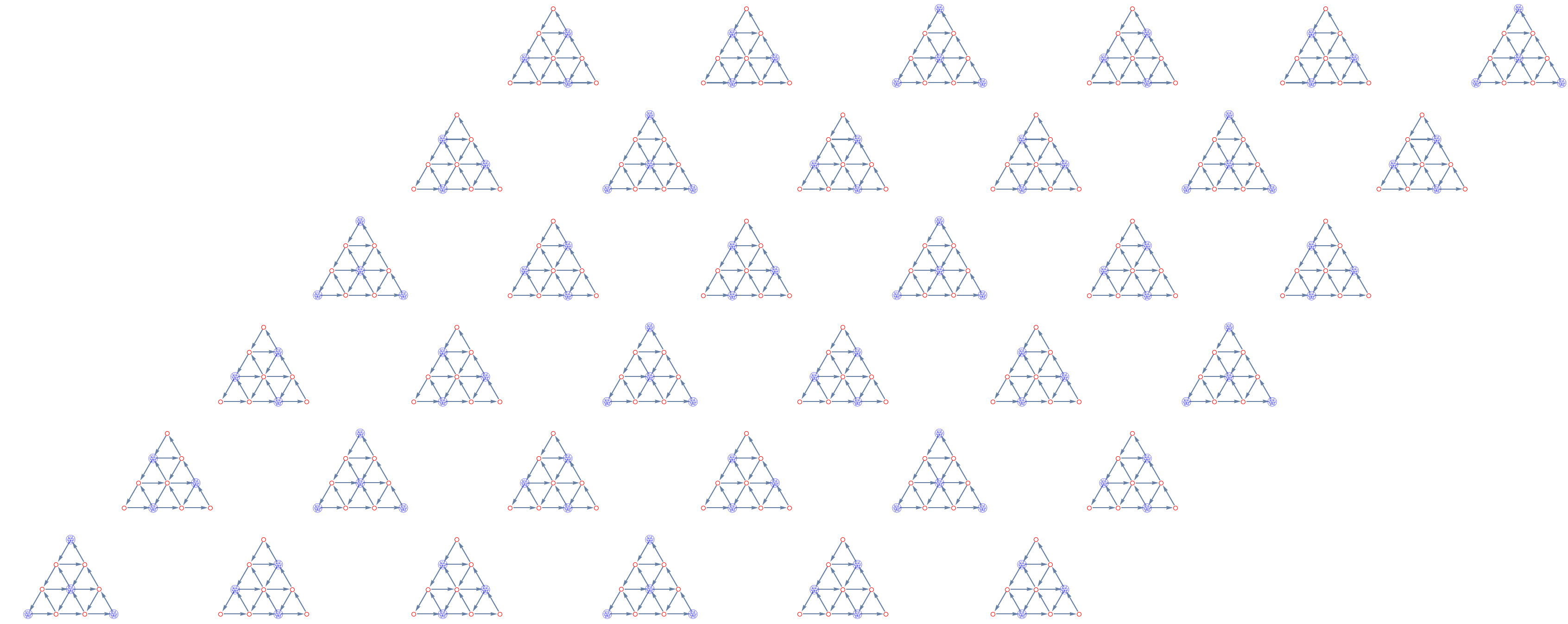}
\caption{\label{A3parallelogramOfRoots} $\SU{3}$ at level $3$ : ${\mathcal A}_3$ parallelogram of roots - not extended.  I should use extended parallelogram.}
\end{figure}
}

\ommit{
\begin{figure}[ht]
\centering
\includegraphics[width=24pc]{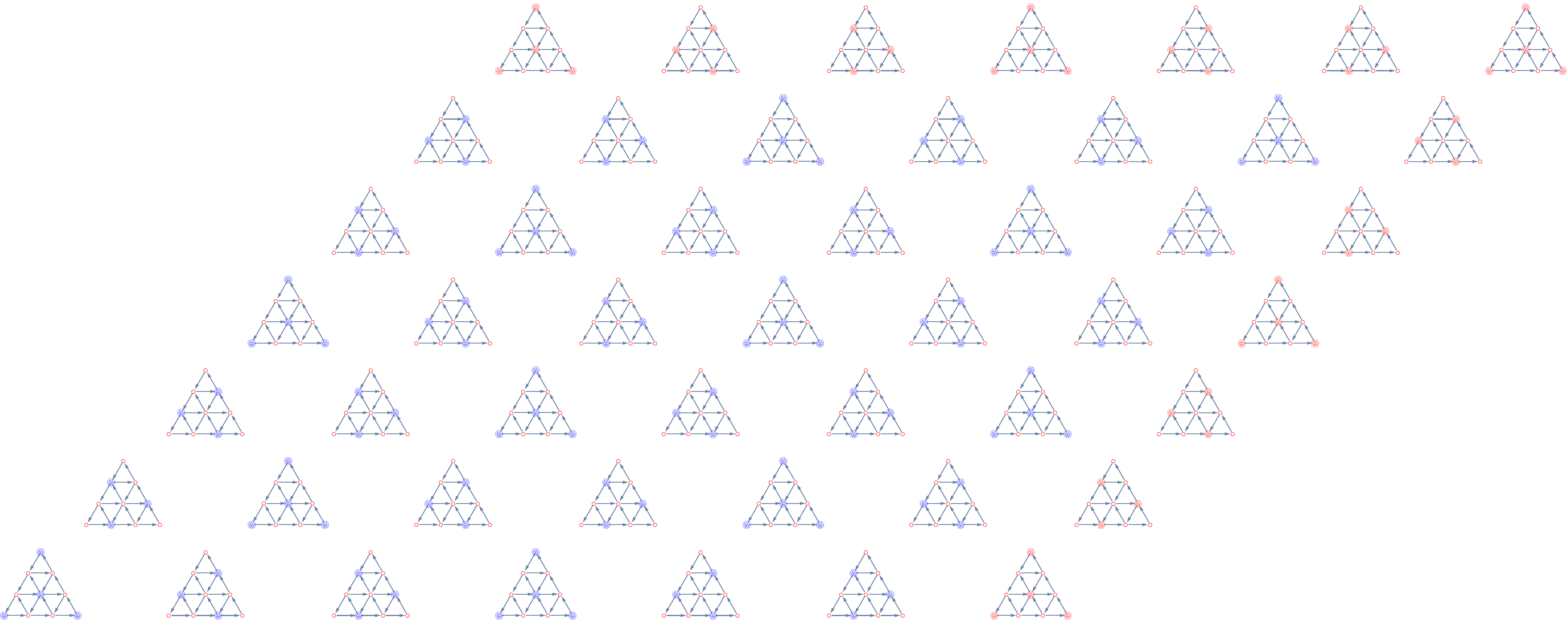}
\caption{\label{ExtendedParallelogramA3} $\SU{3}$ at level $3$ : ${\mathcal A}_3$ extended parallelogram of roots (alone).}
\end{figure}
}

\subsection{Choice of a basis}
\label{choicebasis}
\begin{figure}[ht]
\centering
\includegraphics[width=26pc]{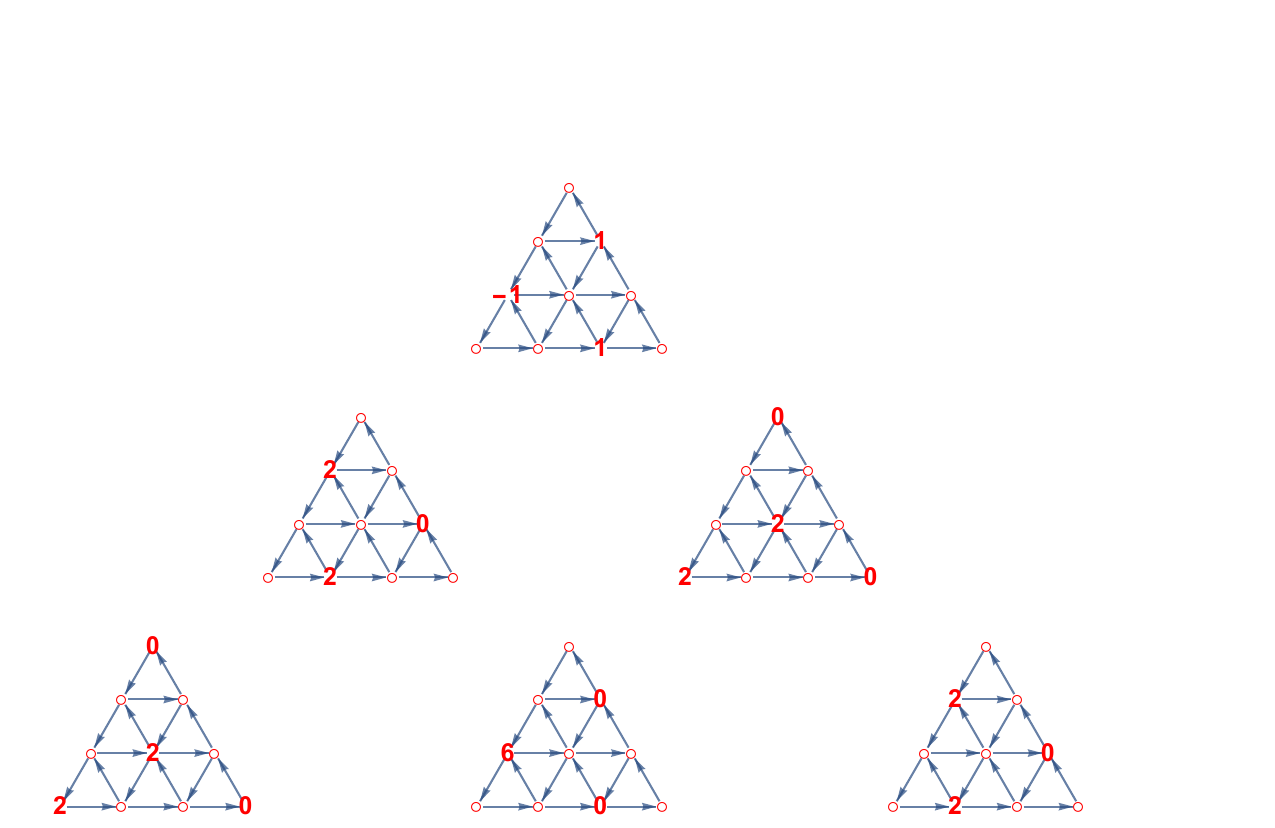}
\caption{\label{A3scalarproductssimplerootsA3} $\SU{3}$ at level $3$: the basis ${\mathcal B}_1$ in the periodicity quiver ($20$ positions marked with integers -- they are located in the bottom left corner of Fig.~\ref{pathLabelsAlpha23_larger}). 
We display the scalar product between the higher root marked $6$ and all other basis elements (these values appear on line $5$ of the Gram matrix for ${\mathcal A}_3$ given in section \ref{latticesproperties}).}
\end{figure}

There are many ways of choosing a basis for a lattice. To every choice is associated a fundamental parallelotope and 
a Gram matrix $A$ (the matrix of inner products in this basis).
In the case of the $\SU{2}$ higher root systems (\ie root systems in the usual sense), one may choose for Gram matrix $A$ the Cartan matrix corresponding to a given Coxeter-Dynkin graph ${\mathcal E}$.
For lattices of higher roots the notion of ``Cartan matrix'' is not available. In what follows we shall present only one Gram matrix, called $A$, since this matrix defines the lattice up to integral equivalence. 
It will be obtained from Eq.~\ref{scalarproductofsu3roots} by choosing a basis that we call ${\mathcal B}_1$.
Its $2 r_{\mathcal E}$ elements (assuming $k>0$) belong to the bottom left corner of the $\SU{3}$ period parallelogram, more precisely,  we choose those higher roots located in the admissible vertices of the six fusion graphs sitting in positions 
 $\{\{0,0\},\{0,1\},\{1,0\},\{1,1\},\{2,0\},\{0,2\}\}$ of the weight lattice; one checks that this indeeds determines a basis which is fully specified once an ordering has been chosen. 
Many other basis choices are of course possible, for instance ${\mathcal B}_2$ or ${\mathcal B}_3$, respectively associated with the admissible vertices belonging to the fusion graphs located in positions
$\{ \{1, 1\}, \{2, 1\}, \{1, 2\}, \{3, 1\}, \{2, 2\}, \{1, 3\}\}$ and $\{\{0,0\},\{1,0\},\{0,1\},\{N-1,N-2\},\{N-2,N-1\},\{N-1,N-1\}\}$, with $N=k+3$, but we shall not display the associated Gram matrices.
The basis ${\mathcal B}_1$ is illustrated in figure~\ref{A3scalarproductssimplerootsA3}.

\subsection{Higher roots and their inner products: summary of the procedure.}
\begin{itemize}
\item Choose a module ${\mathcal E}$ over  ${\mathcal A}_k={\mathcal A}_k(\SU3)$, for instance ${\mathcal A}_k$ itself.
\item From the fundamental fusion matrix $F_{(10)}$ of the chosen module, calculate the other fusion matrices $F_n$,  for instance using the $\SU{3}$ recurrence relation (Eq.~\ref{recF}).
\item Extend the fusion matrices to the weight lattice of $\SU{3}$, using symmetries and periodicity.
\item It is useful to build the periodic essential matrices $\tau_a$, not only the $F_n$, in particular if the chosen module ${\mathcal E}$ is not ${\mathcal A}_k$.
\item Using Eq.~\ref{scalarproductofsu3roots} one can determine a matrix $A^{\mathrm{big}}$ of the $\vert {\mathcal R}\vert \times \vert {\mathcal R}\vert$ scalar products between the $(k+1)(k+2)(k+3)^2 /3$ higher roots (or only between those of ${\mathcal R}^\vee$). The matrix $A^{\mathrm{big}}$ has rank $\mathfrak{r}  = 2 r_{\mathcal E}$.
\item Select a family $(\alpha_i)$ of $\mathfrak{r}$ independent higher roots (\ie choose a basis) and call $A$ the $\mathfrak{r}\times \mathfrak{r}$ restriction of the previous table $A^{\mathrm{big}}$ to the chosen basis.
$A$ will be a Gram matrix for the lattice of higher roots. $A^{\mathrm{big}}$ can however be huge and it is shorter to determine $A$ by calculating only the $\mathfrak{r}\times \mathfrak{r}$ inner products between the basis elements of some chosen basis, for instance ${\mathcal B_1}$, described previously.
\item The choice of $A$ determines a basis $(\alpha_i)$ of higher roots that are such that $< \alpha_i, \alpha_j> = A_{ij}$.\\
Call $K=A^{-1}$ the inverse of $A$ and $(\omega_i)$ the dual basis of $(\alpha_i)$, then  $<\omega_i, \omega_j> = K_{ij}$ and $<~\alpha_i, \omega_j> =~\delta_{ij}$.
The family of vectors $(\omega_i)$ is, by definition, the basis of higher weights associated with the basis of higher roots $(\alpha_i)$.
Linear combinations of the vectors  $(\omega_i)$ with integer coefficients are (integral) higher weights.
Warning: Indices $i,j\ldots$ of $(\alpha_i)$ or of $(\omega_i)$ run from $1$ to ${\mathfrak r}=2 r_{\mathcal E}$ whereas indices $a,b$ of $(\tau_a)$ refer to the irreps of ${\mathcal E}$ and therefore
run only from $1$ to $r_{\mathcal E}$.
\item The last step is to study the lattice of higher roots and its theta function. How this is done is described in the next section.
\end{itemize}

One can a posteriori check that the orthonormal projection of a Dirac measure on the ribbon on the subspace of harmonic functions (higher weights) is indeed a higher root.
This could have been used as a method to determine the latter. 

\section{Theta functions for lattices of higher roots}
\label{resultsforlattices}

\subsection{Lattices and theta functions (reminders)}
\label{ThetaFunctionsTheoryForLattices}
We remind the reader a few results about lattices and their theta functions. This material can be gathered from \cite{Zagier:modularforms}.

Consider a positive definite quadratic form $Q$ which takes integer values on $\ZZ^m$. We can write $Q = \frac{1}{2} x^T \, A \, x$, with $x \in \ZZ^m$ and  $A$ a symmetric $m \times m$ matrix. Integrality of $Q$ implies that $A$ is an even integral matrix (its matrix elements are integers and its diagonal elements are even). Therefore $A$ is a positive definite non singular matrix, and $det(A) >0$. So the inverse $A^{-1}$ exists, as a matrix with rational coefficients. 
The {\sl modular level} of $Q$, or of $A$, is the smallest integer $\ell$ such that  $\ell A^{-1}$ is again an even integral matrix -- this notion differs from the notion of conformal level $k$ used in the previous part of this article. $\Delta = (-1)^m \, det(A)$ is  the {\sl discriminant} of $A$.

Given $Q$, one defines the theta function  $\theta_Q (z) = \sum_{n=0}^\infty  \, p(n) \, q^n$ where\footnote{This parameter $q$ is not related to the root of unity, called ${\mathfrak q}$, that appears in section~\ref{fusioncatbasics}.} $q=exp(2 i \pi z)$  and  $p(n) \in \ZZ_{\geq 0}$  is the number of vectors $x\in Z^m$ that are such that $Q(x)=n$.
The function $\theta_Q$ is always a modular form of weight $m/2$.
In our framework $m$ will always be even (in particular $\Delta = det(A)$) so that we set $m=2s$ with $s$ an integer.

The following theorem  (Hecke-Schoenberg) is known \cite{Zagier:modularforms} and will be used: \\
{\sl
Let $Q : \ZZ^{2s} \mapsto \ZZ$ a positive definite quadratic form, integral, with $m=2s$ variables, of level $\ell$ and discriminant $\Delta$. Then the theta function $\theta_Q$ is a modular form  on the group $\Gamma_0(\ell)$, of weight $s$,  and character $\chi_\Delta$.\\
In plain terms:  $\theta_Q(\frac{az+b}{cz+d}) = \chi_\Delta(a) \, (c z + d)^s \, \theta_Q(z)$ for all  $z\in {\mathfrak H}$ (upper half-plane) and 
$ \left( \begin{smallmatrix} a&b\\ c&d \end{smallmatrix} \right) \in \Gamma_0(\ell)$. Here $\Gamma_0(\ell)$ is the subgroup\footnote{As  $\Gamma_1(\ell) \subset \Gamma_0(\ell)$, one can sometimes use 
modular forms (and bases of spaces of modular forms) twisted by Dirichlet characters on the congruence subgroup $\Gamma_1(\ell)$.}
 of $SL(2,Z)$ defined by the condition $c \equiv 0 \, mod \, \ell$
and $\chi_\Delta$ is the unique Dirichlet character modulo $\ell$ which is such that $\chi_\Delta(p) = {\mathfrak L}(\Delta,p)$ for all odd primes $p$ that do not divide $\ell$, where 
${\mathfrak L}$ denotes the Legendre symbol.}

Notice that $m$, as defined above, is also, in our framework, the dimension of the space ${\mathfrak C}$ of higher roots, which, for $G=\SU{3}$, is equal to $2 r_{\mathcal E}$. 
In that case,  the weight of the (twisted) modular form $\theta_Q$ is therefore equal to $r_{\mathcal E}$, the number of vertices of the fusion diagram, or the number of simple objects in ${\mathcal E} $.

\paragraph{About Dirichlet characters.} Dirichlet characters are particular functions from the integers to the complex numbers that arise as follows:  given a character on the group of invertible elements of the set of integers modulo $p$, one can lift it to a completely multiplicative function on integers relatively prime to $p$ and then extend this function to all integers by defining it to be $0$ on integers having a non-trivial factor in common with $p$. A Dirichlet character with modulus $p$ takes the same value on two integers that agree modulo $p$. The interested reader may consult the abundant literature on the subject but it is enough for us to remember that they are a particular kind of completely multiplicative complex valued functions on the set of integers,  that there are $\phi(p)$ characters modulo $p$, where $\phi$ is the Euler function, and that they are tabulated in many places ---there is even a command DirichletCharacter[p,j,n] in Mathematica \cite{Mathematica} that gives the Dirichlet character with modulus $p$ and index $j$ as a function of $n$ (the index $j$ running from $1$ to $\phi(p)$).

\subsection{Lattice properties: tables}
From the fusion matrices associated with the different quantum subgroups one determines Gram matrices for the associated  lattices of higher roots.
For the first few members of the series considered in this paper, and for the exceptional cases, some information, encoded in the Gram matrices, is summarized in the following table\footnote{In the table, the subscript $x$ of $\Gamma_x$ may be $0$ or $1$ (see footnote 10) and the entry  \text {kiss}
gives the smallest term of the theta series, its coefficient being the kissing number.}

\bigskip

\begin{tabular}{L | L  L  L L L L L L L L L}
{\mathcal E} & k & N & r_{\mathcal E} & \mathfrak{r} & \vert {\mathcal R} \vert & \text {kiss} & \Delta & \ell & \varphi_{Euler}(\ell) & \text {dim}(M_{ r_{\mathcal E} }(\Gamma_x(\ell),\chi)) & \vert aut \vert \\
\hline
{\mathcal A}_1& 1 & 4 & 3 & 6 & 32 & 32 \, q^6 & 4^6 & 16 & 8 & 7 & 6!/2 \times 64\\
{\mathcal A}_2& 2 & 5 & 6 & 12 & 100 & 100 \, q^6 & 5^9 & 25 & 20 & 16 & 1200 \\
{\mathcal A}_3& 3 & 6 & 10 & 20 & 240 & 240 \, q^6 & 6^{12} & 18 & 6 & 31 & 864 \\
{\mathcal A}_4& 4 & 7 & 15 & 30 & 490 & 490 \, q^6 & 7^{15} & 49 & 42  & 70 & 3528\\
\hline
{\mathcal D}_3& 3 & 6 & 6 & 12 & 144 & 36 \, q^4 & 3^{12} & 9 & 6 & 7 & 6912 \\
{\mathcal D}_6& 6 & 9 & 12 & 24 & 648 & 162 \, q^4 & 3^{18} & 27 & 18 & 36 & 2^6 3^{11}  \\
\hline
{\mathcal E}_5& 5 & 8 & 12 & 24 & 512 & 512 \, q^6 & 2^{30} & 16 & 8 & 25  & 2^{11} 3 \\
{\mathcal E}_9& 9 & 12 & 12 & 24 & 1152 & 756 \, q^4 & 2^{24} & 16 & 8 & 13  &2^{10} 3^4 5 \\
{\mathcal E}_{21}&21 & 24 & 24 & 48 & 9216 & 144 \, q^4 & 3^{12} & 3 & 2 & 9 &  2^{12} 3^2
\end{tabular}

We give below the coefficients of the associated theta series in the variable $q^2$ (first term is $(q^2)^0=1$)), for instance,
$\theta({\mathcal A}_2)=1+ 100 (q^2)^3 +450 (q^2)^4 + \ldots = 1+100 q^6 + 450 q^8 + 960 q^{10} +2800 q^{12}+  \ldots$.

\small
\begin{equation*}
\begin{split}
{\mathcal A}_0 (\text{rescaled}): \; &
1, 6, 0, 6, 6, 0, 0, 12, 0, 6, 0, 0, 6, 12, 0, 0, 6, 0, 0, 12, 0, 12, 0, 0, 0, 6, 0, 6, 12, 0, 0, 12, 0, 0, 0, 0, 6, 12, 0, 12,\\
& 0, 0, 0, 12, 0, 0, 0, 0, 6, 18, 0, 0, 12, 0, 0, 0, 0, 12, 0, 0, 0, 12, 0, 12, 6, 0, 0, 12, 0, 0, 0, 0, 0, 12, 0, 6, 12, 0, 0, 12, 0,\ldots\\
{} & {} \\
{\mathcal A}_1: \; &
1, 0,0,32,60,0,0,192,252,0,0,480,544,0,0,832,1020,0,0,1440,1560,0,0, 2112, 2080, 0, 0, 2624, 3264,\\
& 0, 0, 3840, 4092, 0, 0, 4992, 4380, 0, 0, 5440, 6552, 0, 0, 7392, 8160, 0, 0, 8832, 8224,\ldots \\
{} & {} \\
{\mathcal A}_2: \; &
1,0,0,100, 450, 960, 2800, 6600, 12300, 22400, 30690, 63000, 93150, 144000, 203100, 236080, 392850, \\
& 550800, 708350, 961800, 972780, 1581600, 1937250, 2495400, 2977400, 3063360, 4469400, 5547700, \\
& 6477600, 7963200, 7344920, 11094000, 12627000, 15127200, 17091900, 16459440, 22670850, 26899200,\ldots\\
{} & {} \\
{\mathcal A}_3: \; &
1,0,0,240,1782,9072,59328,216432,810000,2059152,6080832,12349584,31045596,57036960,122715648, \\
& 204193872,418822650,622067040,1193611392,1734272208,3043596384,4217152080,7354100160,9446435136, \\
& 15901091892,20507712192,32268036096,40493364288,64454759856,76079125584,118436670720,142127536464,\ldots\\
{} & {} \\
{\mathcal A}_4: \; &
1, 0, 0, 490, 4998, 45864, 464422, 3429426, 21668094, 111678742, 492567012, 1876801038,\\
& 6352945942, 19484903508, 54935857326, 144330551050,\ldots\\
{} & {} \\
{\mathcal A}_5: \; &
1, 0, 0, 896, 11856, 154368, 2331648, 27065088, 281311128, \ldots\\
{} & {} \\
{\mathcal A}_6: \; &
1,0,0,1512, 24300 ,425736, 8530758, \ldots\\
{} & {} \\
{} & {} \\
{\mathcal D}_3: \; &
1, 0, 36, 144, 486, 2880, 5724, 7776, 31068, 40320, 47628, \ldots\\
{} & {} \\
{\mathcal D}_6: \; &
1, 0, 162, 2322, 35478, 273942, 1771326, 9680148, 40813632, 150043014, 484705782\ldots\\
{} & {} \\
{} & {} \\
{\mathcal E}_5: \; &
1, 0, 0, 512, 11232, 145920, 1055616, 5618688, 25330128, 89127936, 295067136,
810542592, 2185379968, 5109275136,
\ldots\\
{} & {} \\
{\mathcal E}_9: \; &
1, 0, 756, 5760, 98928, 1092096, 8435760, 45142272, 202712400, 715373568, 2350118808, 
6501914496, 17469036096,
\ldots\\
{} & {} \\
{\mathcal E}_{21}: \; &
1, 0, 144, 64512, 54181224,\ldots\\
\end{split}
 \end{equation*}
 \normalsize
\vfill
\eject

\subsection{Properties of the lattices}
\label{latticesproperties}
\paragraph{Case ${\mathcal A}_0$.}
This lattice coincides with the (rescaled) usual root lattice of $\SU{3}$ also known as the planar hexagonal lattice, and ${\mathcal R}^\vee$ with the set of positive roots.
Although its theta function can be found in many textbooks\footnote{Its expression in terms of the elliptic theta function $\vartheta_3$ reads:
 \begin{eqnarray*}
 \theta(z) &=& \frac{\vartheta _3(0,q){}^3+\vartheta _3\left(\frac{\pi }{3},q\right){}^3+\vartheta
   _3\left(\frac{2 \pi }{3},q\right){}^3}{3 \vartheta _3\left(0,q^3\right)}
    \end{eqnarray*}}, for instance in \cite{ConwaySloane},
  it is instructive to obtain it by using the theorem recalled in the previous section.
  The Gram matrix, in the basis of higher roots obtained from Eq.~(\ref{scalarproductofsu3roots}), is
 $\left(\begin{smallmatrix}6&-3\\-3&6\end{smallmatrix}\right)$, \ie three times the Cartan matrix of $\SU{3}$.
 Its theta function is therefore a modular form on the group $\Gamma_0(3)$, of weight $s=1$ twisted by the (unique in this case) non-trivial Dirichlet character modulo $3$; this vector space of modular forms is of dimension $1$, hence $\theta$ can be identified with its generator.

\paragraph{Case ${\mathcal A}_1$.}
The Gram matrix $A$, with the basis choice ${\mathcal B}_1$, and its inverse $K$ are given below.
\be
\begin{array}{ccc}
A=
\left(
\begin{array}{cccccc}
 6 & 2 & 2 & -2 & -2 & -2 \\
 2 & 6 & 2 & 2 & -2 & 2 \\
 2 & 2 & 6 & 2 & 2 & -2 \\
 -2 & 2 & 2 & 6 & 2 & 2 \\
 -2 & -2 & 2 & 2 & 6 & -2 \\
 -2 & 2 & -2 & 2 & -2 & 6 \\
\end{array}
\right)
& &
K= \frac{1}{8} \;
\left(
\begin{array}{cccccc}
 3 & -1 & -1 & 1 & 1 & 1 \\
 -1 & 3 & -1 & -1 & 1 & -1 \\
 -1 & -1 & 3 & -1 & -1 & 1 \\
 1 & -1 & -1 & 3 & -1 & -1 \\
 1 & 1 & -1 & -1 & 3 & 1 \\
 1 & -1 & 1 & -1 & 1 & 3 \\
\end{array}
\right)
\end{array}
\ee

The $16$ elements of ${\mathcal R}^\vee$ may be called ``positive higher roots'' (their opposites, the elements of $-{\mathcal R}^\vee$, being ``negative''), they can be expanded on the root basis
${\mathcal B}_1=\{{\alpha_i}\}, i=1\ldots 6$ as follows:
{\footnotesize
\begin{equation*}
\begin{split}
\alpha _ 1 - \alpha _ 2 + \alpha _ 6,\,&
 -\alpha _ 2 + \alpha _ 3 - \alpha _ 5,\, -\alpha _ 1 + \alpha _ 3 - \alpha _ 4,\, -\alpha _ 4 + \alpha _ 5 + \alpha _ 6,\, \alpha_ 6,\, -\alpha _ 2 + \alpha _ 4 - \alpha _ 5,\, -\alpha _ 1 - \alpha _ 5 - \alpha _ 6,  \\
  -\alpha _ 1 + \alpha _ 2 - \alpha _ 4,&
   \, \alpha _ 2,\, \alpha _ 4,\, -\alpha _ 3 + \alpha_ 4 - \alpha _ 6,\, \alpha _ 2 - \alpha _ 3 - \alpha _ 6, \, \alpha _ 1,\, \alpha _ 3,\, \alpha _ 5,\, \alpha _ 1 - \alpha _ 3 + \alpha _ 5
\end{split}
\end{equation*}}
With the same ordering,  the family of their mutual inner products builds a  $16\times 16$ matrix  $A^{\mathrm{big}}$,
which is of rank $6$, as expected.

The lattice  is even, with minimal norm~$6$, and all the coefficients of the above Gram matrix are even; if we rescale it, setting  $B = A/2$, the vectors of minimal norm have then norm $3$. 
The determinant of $A$, or ``connection index'', is also the order of the dual quotient, an abelian group isomorphic with $Z_2\times (Z_4)^{ \times 4} \times Z_8$. The  lattice defined by $A$ is obviously not self-dual. If we use its rescaled version, the connection index becomes $64$ and the dual quotient is then isomorphic with $(Z_2)^{ \times 4} \times Z_4$; elements of the (rescaled) dual belong to one and only one congruence class, an element of the dual quotient, and are  therefore be classified by 5-uplets $(c_{2_1},c_{2_2},c_{2_3},c_{2_4},c_4)$, with $c_{2_i} \in \{0,1\}$ and $c_4  \in \{0, 1, 2, 3\}$.

\smallskip

We find the theta function of this lattice by applying the Hecke-Schoenberg theorem.
From the Gram matrix one finds that the discriminant  is $4^{6}$ and that the  (modular) level of the quadratic form is $16$. 
The odd primes not dividing $16$ are $3, 5, 7, 11, 13$ and their Legendre symbols are all equal to $1$.
From the $8 \times 16$ table of Dirichlet characters of modulus $16$ over the cyclotomic field of order $\varphi_{Euler}(16)=8$
restricted to odd primes not dividing the level,
one selects the unique character whose values 
coincide with the list obtained for the Legendre symbols.
The space of modular forms on $\Gamma_ 1(16)$ of weight $3$,  twisted by this Dirichlet character,
namely the Kronecker  character -4,  has dimension $7$. It is spanned by the following forms (we set $q_2 = q^2$): 
{\small
\begin{eqnarray*}
  b_1&=&  1 + 12 \, q_2^8 + 64 \, q_2^{12} + 60 \, q_2^{16} + O(q_2^{24}), \qquad
  b_2 =  \, q_2 + 21 \, q_2^9 + 40 \, q_2^{13} + 30 \, q_2^{17} + 72 \, q_2^{21} + O(q_2^{24}),\\
  b_3&=&  \, q_2^2 + 26 \, q_2^{10} + 73 \, q_2^{18} + O(q_2^{24}), \qquad
  b_4 =  \, q_2^3 + 6 \, q_2^7 + 15 \, q_2^{11} + 26 \, q_2^{15} + 45 \, q_2^{19} + 66 \, q_2^{23} + O(q_2^{24}),\\
  b_5&=&  \, q_2^4 + 4 \, q_2^8 + 8 \, q_2^{12} + 16 \, q_2^{16} + 26 \, q_2^{20} + O(q_2^{24}),\quad
  b_6 =   \, q_2^5 + 2 \, q_2^9 + 5 \, q_2^{13} + 10 \, q_2^{17} + 12 \, q_2^{21} + O(q_2^{24}),\\
  b_7&=&  \, q_2^6 + 6 \, q_2^{14} + 15 \, q_2^{22} + O(q_2^{24})
\end{eqnarray*}
}
An explicit determination of the vectors (and their norms) belonging to the first shells shows that the theta function starts as $1+32 q_2^3+60 q_2^4 + O(q^{14})$.
The components of this modular form on the previous basis are therefore ${1, 0, 0, 32, 60, 0, 0}$. 
In other words, $$\theta = b_1 + 32 \, b_4 + 60 \, b_5$$
Using a computer package, one can quickly obtain the $q$-expansion of the functions $b_n$ to very large orders. 
 Here is a Magma \cite{Magma} program  that returns its series expansion up to order 48 in $q_2$ and uses the above ideas: 
 {\footnotesize
\begin{verbatim}
H := DirichletGroup(16,CyclotomicField(EulerPhi(16))); 
chars := Elements(H); eps := chars[2];
M := ModularForms([eps],3); order:=48;
PowerSeries(M![1,0,0,32,60,0,0],order);
\end{verbatim}}

One finds that the automorphism group $aut$ of this lattice is of order $23040$ and that it is  isomorphic with the semi-direct product of $A_6$ (the alternated group of order $6!/2= 360$) times an abelian group of order  $64$, actually with
$((C_2)^{\times 5} \rtimes A_6)\rtimes C_2$. Orbits of the basis vectors under the $aut$ action coincide and contain the $32$ higher roots (the $16$ positive and the $16$ negative ones). 
From the given Gram matrix and using for example Magma, one finds that the Voronoi polytope has $92$ $5$-dimensional facets, $4896$ edges and $588$ vertices. 
\ommit{
The group $aut$ is generated by the following matrices 
{\footnotesize
$$
\begin{array}{ccc}
\left(
\begin{array}{cccccc}
 0 & 0 & 0 & -1 & 1 & -1 \\
 -1 & 0 & 0 & -1 & 0 & -1 \\
 -1 & 0 & 1 & 0 & 0 & -1 \\
 0 & 0 & 0 & 0 & -1 & 0 \\
 -1 & 1 & 0 & 0 & -1 & 0 \\
 0 & 0 & 0 & 1 & 0 & 0
\end{array}
\right)
, &
\left(
\begin{array}{cccccc}
 1 & -1 & 0 & 0 & 1 & 0 \\
 0 & -1 & 0 & 0 & 1 & -1 \\
 0 & -1 & 0 & 0 & 0 & 0 \\
 -1 & 0 & 1 & 0 & 0 & -1 \\
 -1 & 0 & 0 & -1 & 0 & -1 \\
 0 & 0 & 0 & 0 & -1 & 0
\end{array}
\right)
, &
\left(
\begin{array}{cccccc}
 1 & 0 & 0 & 0 & 0 & 0 \\
 0 & 1 & 0 & 0 & 0 & 0 \\
 0 & 0 & 0 & -1 & 1 & -1 \\
 0 & 0 & 0 & 1 & 0 & 0 \\
 0 & 0 & 0 & 0 & 1 & 0 \\
 0 & 0 & -1 & -1 & 1 & 0
\end{array}
\right)
\end{array}
$$
}
}

{\sl Other avatars of this lattice:}
The obtained theta series starts as the theta series of a  (scaled version of) the shifted  $D_6$ lattice, called $D_6^{+} = D_6 \cup ([1] + D_6)$, see \cite{ConwaySloane}.
This coincidence (already noticed in \cite{Ocneanu:Bariloche}) is not sufficient to allow an identification with $D_6^{+}$, but it is so because one can choose the same Gram matrix for both lattices.
Using this identification, we can re-write the theta series in terms of elliptic theta functions as follows:
$$\frac{1}{2} \left(\vartheta _2\left(0,q^4\right){}^6+\vartheta _3\left(0,q^4\right){}^6+\vartheta _4\left(0,q^4\right){}^6\right)$$
This an alternative to the expression of $\theta$  previously given in terms of appropriate modular forms.
It is known \cite{ConwaySloane} that the $D_n^{+}$ packing is a lattice packing if and only if $n$ is even. In particular it is so for $n=6$. 
The fact that $D_n^{+}$ is not a lattice for $n$ odd excludes a possible systematic identification with lattices of type ${\mathcal A}_k(\SU3)$ when $k>1$.
Here are a few others avatars of $D_6^{+}$ :\\
- The generalized laminated lattice $\Lambda_6[3]$ with minimal norm $3$, see \cite{PleskenPohst} (the authors study the family $\Lambda_n[3]$ and provide enough information to allow one to recover the Gram matrix associated with $\Lambda_6[3]$ and show that it coincides with the matrix $A$ given before).\\
- The lattice ${\mathcal L}_4$ generated by cuts of the complete graph on a set of $4$ vertices, see  \cite{DezaGrishukhin} (the authors are interested in the Delaunay polytopes for the lattices ${\mathcal L}_n$; when $n=4$, this lattice is isomorphic with $D_6^{+}$, actually, the precise relation is ${\mathcal L}_4 = {\sqrt 2} \, D_6^{+}$).\\
Identification of lattices defined by ${\mathcal A}_k(\SU3)$ with other members of the above families fails.

\paragraph{The case ${\mathcal A}_2$.}
We obtain the following Gram matrix.
\be
A = \left(
\begin{array}{cccccccccccc}
 6 & 0 & 2 & 0 & 2 & 0 & -2 & 1 & -2 & 2 & -2 & 2 \\
 0 & 6 & 2 & 2 & 2 & 2 & 1 & -1 & 0 & -2 & 0 & -2 \\
 2 & 2 & 6 & 0 & 2 & 2 & 2 & 2 & -1 & 1 & 2 & 2 \\
 0 & 2 & 0 & 6 & 2 & 0 & 0 & 2 & 1 & -2 & 2 & 0 \\
 2 & 2 & 2 & 2 & 6 & 0 & 2 & 2 & 2 & 2 & -1 & 1 \\
 0 & 2 & 2 & 0 & 0 & 6 & 0 & 2 & 2 & 0 & 1 & -2 \\
 -2 & 1 & 2 & 0 & 2 & 0 & 6 & 0 & 2 & 0 & 2 & 0 \\
 1 & -1 & 2 & 2 & 2 & 2 & 0 & 6 & 2 & 2 & 2 & 2 \\
 -2 & 0 & -1 & 1 & 2 & 2 & 2 & 2 & 6 & 0 & 0 & -2 \\
 2 & -2 & 1 & -2 & 2 & 0 & 0 & 2 & 0 & 6 & -2 & 2 \\
 -2 & 0 & 2 & 2 & -1 & 1 & 2 & 2 & 0 & -2 & 6 & 0 \\
 2 & -2 & 2 & 0 & 1 & -2 & 0 & 2 & -2 & 2 & 0 & 6 \\
\end{array}
\right)
\label{GramA2}
\ee
The discriminant is readily calculated:  $\Delta = 5^9$. The modular level is $\ell = N^2 = 25$. 
Applying the Hecke-Schoenberg theorem leads to the following result:
the theta function of this lattice is of weight $6$, modular level $\ell = 5^2 = 25$ (the square of the altitude) and Dirichlet character $\chi(11)$ for the characters modulo 25 on a cyclotomic field of order 20. It is the only character (namely the Kronecker character 5), the eleventh on a collection of $20 = \Phi_{Euler}(25)$) that coincides with the value of the Legendre symbol ${\mathfrak L}(\Delta,p)$ for all odd primes $p$ that do not divide $25$. This space of modular forms has dimension $16$. The theta function, in the variable $q_2=q^2$, is therefore fully determined by its $16$ first Fourier coefficients (the first being $1$). The coefficients of $q_2^a$ with $a>15$ are then predicted. 
Here is the Magma code calculating the first $48\times 2$ coefficients:
{\footnotesize
\begin{verbatim}
H := DirichletGroup(25,CyclotomicField(EulerPhi(25)));
chars := Elements(H); eps := chars[11];
M := ModularForms([eps],6); order:=48;
PowerSeries(M![1, 0, 0,100, 450, 960, 2800, 6600, 12300, 22400, 30690, 63000, 
93150, 144000, 203100, 236080],order);
\end{verbatim}}
The first Fourier coefficients have to be computed by a brute force approach that relies, ultimately, on the obtained Gram matrix.
The automorphism group of this lattice is of order $1200$ and its structure, in terms of direct and semi-direct products, is $C_2 \times ((((C_5 \times C_5) \rtimes C_4) \rtimes C_3) \rtimes C_2)$.
Orbits of the basis vectors under the $aut$ action coincide and contain the $100$ higher roots (the $50$ positive and the $50$ negative ones). 
Their stabilizers are conjugated in $aut$, and are isomorphic with the group $D_{12}$ (which is itself isomorphic with $S_3 \times C_2$).
We  only mention that the Voronoi polytope has $5410$ $11$-dimensional facets.

\paragraph{The case ${\mathcal A}_3$.}
A Gram matrix is given below. 
We only mention that the discriminant is $\Delta = 6^{12}$, the modular level is $\ell = 18$,  
the rank of the lattice $L_3 $ is $\mathfrak{r} = 2\times\, r_{{\mathcal A}_3} = 20$,  the period is a rhombus $6\times 6$ and  $\vert {\mathcal R}\vert = 240$.
The theta function belongs to  a space of modular forms on $\Gamma_0(18)$, of weight $10$, twisted by an appropriate character of modulus $18$ on a cyclotomic field of order $6 = \Phi_{Euler}(18)$. The corresponding space of modular forms has dimension $31$ and the theta function of the lattice is fully determined by its first Fourier coefficients. The automorphism group has order $864$ and the Voronoi polytope has $539214$ $19$-dimensional facets.
{\scriptsize
\be 
A=
\left(
\begin{array}{cccccccccccccccccccc}
 6 & 0 & 0 & 0 & 2 & 0 & 0 & 2 & 0 & 0 & -2 & 1 & 0 & 0 & -2 & 2 & 0 & -2 & 2 & 0 \\
 0 & 6 & 0 & 0 & 2 & 2 & 2 & 2 & 2 & 2 & 1 & 0 & 1 & 1 & 0 & 0 & 0 & 0 & 0 & 0 \\
 0 & 0 & 6 & 0 & 0 & 0 & 2 & 0 & 2 & 0 & 0 & 1 & -2 & 0 & 2 & 0 & -2 & 0 & -2 & 2 \\
 0 & 0 & 0 & 6 & 0 & 2 & 0 & 0 & 0 & 2 & 0 & 1 & 0 & -2 & 0 & -2 & 2 & 2 & 0 & -2 \\
 2 & 2 & 0 & 0 & 6 & 0 & 0 & 2 & 2 & 0 & 2 & 2 & 0 & 0 & -1 & 1 & 1 & 2 & 2 & 0 \\
 0 & 2 & 0 & 2 & 0 & 6 & 0 & 2 & 0 & 2 & 0 & 2 & 0 & 2 & 1 & -1 & 1 & 2 & 0 & 2 \\
 0 & 2 & 2 & 0 & 0 & 0 & 6 & 0 & 2 & 2 & 0 & 2 & 2 & 0 & 1 & 1 & -1 & 0 & 2 & 2 \\
 2 & 2 & 0 & 0 & 2 & 2 & 0 & 6 & 0 & 0 & 2 & 2 & 0 & 0 & 2 & 2 & 0 & -1 & 1 & 1 \\
 0 & 2 & 2 & 0 & 2 & 0 & 2 & 0 & 6 & 0 & 0 & 2 & 2 & 0 & 2 & 0 & 2 & 1 & -1 & 1 \\
 0 & 2 & 0 & 2 & 0 & 2 & 2 & 0 & 0 & 6 & 0 & 2 & 0 & 2 & 0 & 2 & 2 & 1 & 1 & -1 \\
 -2 & 1 & 0 & 0 & 2 & 0 & 0 & 2 & 0 & 0 & 6 & 0 & 0 & 0 & 2 & 0 & 0 & 2 & 0 & 0 \\
 1 & 0 & 1 & 1 & 2 & 2 & 2 & 2 & 2 & 2 & 0 & 6 & 0 & 0 & 2 & 2 & 2 & 2 & 2 & 2 \\
 0 & 1 & -2 & 0 & 0 & 0 & 2 & 0 & 2 & 0 & 0 & 0 & 6 & 0 & 0 & 0 & 2 & 0 & 2 & 0 \\
 0 & 1 & 0 & -2 & 0 & 2 & 0 & 0 & 0 & 2 & 0 & 0 & 0 & 6 & 0 & 2 & 0 & 0 & 0 & 2 \\
 -2 & 0 & 2 & 0 & -1 & 1 & 1 & 2 & 2 & 0 & 2 & 2 & 0 & 0 & 6 & 0 & 0 & 0 & -2 & 2 \\
 2 & 0 & 0 & -2 & 1 & -1 & 1 & 2 & 0 & 2 & 0 & 2 & 0 & 2 & 0 & 6 & 0 & -2 & 2 & 0 \\
 0 & 0 & -2 & 2 & 1 & 1 & -1 & 0 & 2 & 2 & 0 & 2 & 2 & 0 & 0 & 0 & 6 & 2 & 0 & -2 \\
 -2 & 0 & 0 & 2 & 2 & 2 & 0 & -1 & 1 & 1 & 2 & 2 & 0 & 0 & 0 & -2 & 2 & 6 & 0 & 0 \\
 2 & 0 & -2 & 0 & 2 & 0 & 2 & 1 & -1 & 1 & 0 & 2 & 2 & 0 & -2 & 2 & 0 & 0 & 6 & 0 \\
 0 & 0 & 2 & -2 & 0 & 2 & 2 & 1 & 1 & -1 & 0 & 2 & 0 & 2 & 2 & 0 & -2 & 0 & 0 & 6 \\
\end{array}
\right)
\label{GramA3}
\ee
}
\paragraph{Cases ${\mathcal D}_3, {\mathcal D}_6,\ldots, {\mathcal E}_5, {\mathcal E}_9, {\mathcal E}_{21}$.} The procedure should be clear by now and we just refer to the previously given tables.
Explicit Gram matrices, in particular for the three exceptional cases, can be found in one appendix of \cite{RC:hyperroots}.

\subsection{Remarks}

\paragraph{About the vectors of smallest norm.}  For the lattices associated with ${\mathcal A}_k(\SU{3})$ that we considered explicitly, 
 the lattice vectors of shortest length are precisely the higher roots ($100$ of them for $L_3$, for instance), the kissing number of those lattices are then given by the number of higher roots. 
 As it is well known, this property holds for all usual root lattices, \ie higher root lattices of the $\SU{2}$ family.
However this property does not always hold for those lattices associated with modules  of the $\SU{3}$ family that are not of type ${\mathcal A}_k$ although it holds for ${\mathcal E}_5$.
 In the case ${\mathcal D}_{3}$ the first shell is made of vectors of norm $4$, so they are not higher roots, and the only vectors of the lattice that belong to the second shell, of norm $6$, are precisely the higher roots. 
 In the case ${\mathcal D}_{6}$, like for ${\mathcal E}_9$ or ${\mathcal E}_{21}$  the first shell is made of vectors of norm $4$ (which are not higher roots), but the second shell, of norm $6$, contains not only the higher roots themselves, but other vectors as well. In all cases the vectors of smallest norm can of course be expanded on a chosen basis of higher roots.

\paragraph{About the determination of $\theta({\mathcal A}_k)$, for general $k$.}
The theta function, as a modular form twisted by a character, can, in principle, be obtained by following the method explained in the previous sections and illustrated in the case of the first few members of the ${\mathcal A}_k$ series.
In this respect we observed that the (quadratic form) level of is often equal to  $\ell = (k+3)^2$ but it is not always so.
The discriminant is $(k + 3)^{3 (k + 1)}$,
the weight is $r_{{\mathcal A}_k} = (k+1)(k+2)/2$, the quadratic form level is readily obtained from the Gram matrix, and the determination of the appropriate character requires a discussion relying on the arithmetic properties of the discriminant and of the level. 
However, the first coefficients of the Fourier series expansion have to be found, and the number of needed coefficients depends on the properties of an appropriate space of modular forms. The determination of the needed coefficients is done by brute force, namely by computing the norm of the vectors belonging to the first shells, using the Gram matrix as an input.  Moreover, the explicit determination of a Gram matrix becomes a non-trivial exercise when $k$ is large.
The present method may therefore become rapidly intractable if we increase $k$ too much. 
\ommit{
Admittedly it would be nice to have a general formula, like the one that we have for the root lattices of type $A_{n-1}$, that would be valid for all $k$'s, and would express the theta function of $L_k$ in terms of known functions (for instance elliptic theta's). This was not done but we hope that our results will trigger new developments in that direction.}
\paragraph{Construction of lattices of higher roots for quantum subroups of type  $(G,k)$, for other Lie groups.} The calculations described in the present paper, leading to explicit theta series for lattices of higher roots associated with quantum subgroups of type $(\SU3,k)$ could be generalized, without much ado, to the other Lie groups of rank $2$, namely $B_2$ and $G_2$, since the fusion matrices of their quantum modules and quantum subgroups at level $k$ are available, see \cite{CTR: B2G2}, see also \cite{RCsiteWebFusionGraphs} (the exceptional cases given there are obtained by conformal embeddings but it is believed that, at least for these Lie groups,  there are no others). The list of quantum modules and subgroups of type $(\SU4,k)$ is also known --- see \cite{Ocneanu:Bariloche} and \cite{CoquereauxSchieberSU4}.


\begin{thebibliography}{99}
\footnotesize 

\bibitem{CIZ}  Cappelli A.,  Itzykson  C.  and   Zuber J. -B.,  The ADE classification of minimal and $A_{1}^{(1)}$ conformal invariant theories, {\it Commun.  Math.  Phys.},  {13},  pp 1--26, (1987).  
\bibitem{ConwaySloane} Conway J. and Sloane N.J.A., Sphere Packings, Lattices and Groups (3rd ed.), Springer, (1999).
 \bibitem{RC:hyperroots} Coquereaux R., {Theta functions for lattices of $\SU{3}$ hyper-roots},  arXiv:1708.00560,  Experimental Mathematics, 29:2, 137-162, (2020, published online: 02 Apr 2018), DOI: 10.1080/10586458.2018.1446062 
\bibitem{RC:periodicquivers} Coquereaux R., {Quantum McKay correspondence and global dimensions for fusion and module-categories associated with Lie groups}, {\it Jour. of Algebra}, {398}, pp 258-283, (2014).
\bibitem{CordobaRC} Coquereaux R., {Global dimensions for Lie groups at level $k$ and their conformally exceptional quantum subgroups}, {\it Revista de la Union Matematica Argentina}, {51 No 2}, pp 17-42 (2010).
\bibitem{CoquereauxSchieberJMP} Coquereaux R. and Schieber G., {Orders and dimensions for sl2 or sl3 module-categories and boundary conformal field theories on a torus},   {\it J. of Mathematical Physics} {\bf 48} (2007) 043511;
\url{http://arxiv.org/abs/math-ph/0610073}
\bibitem{CoquereauxSU3Maroc}Coquereaux R., Hammaoui D., Schieber G., Tahri  E.H., {Comments about quantum symmetries of $\SU{3}$ graphs},  {\it Journal of Geometry and Physics} 57 pp 269-292 (2006).
\bibitem{RCsiteWebFusionGraphs} Coquereaux R., Fusion graphs, \url{http://www.cpt.univ-mrs.fr/~coque/quantumfusion/FusionGraphs.html}
\bibitem{CIG:TriangularCells}Coquereaux R., Isasi E., Schieber G., {Notes on TQFT wire models and coherence equations for $\SU{3}$ triangular cells}, arXiv :1007.0721. {\it Symmetry, Integrability and Geometry : Methods and Applications}, SIGMA 6 (2010), 099, 44 pp.
\bibitem{CTR: B2G2} Coquereaux R., Tahri E.H., Rais R., { Exceptional quantum subgroups for the rank two Lie algebras B2 and G2}, arXiv:1001.5416. {\it Journal of Mathematical Physics}, Vol.51, Issue 9 (2010).
\bibitem{CoquereauxZuberNuclPhys} Coquereaux R and  Zuber J.-B.,  {On some properties of $\SU{3}$ Fusion Coefficients.},
Contribution to Mathematical Foundations of Quantum Field Theory, special issue in memory of Raymond Stora, 33 pp., {\it Nucl. Phys. B.}, DOI: 10.1016/j.nuclphysb.2016.05.029, (2016).
 \bibitem{CoquereauxSchieberSU4}  Coquereaux R. and Schieber G., From conformal embeddings to quantum symmetries: an exceptional $\SU{4}$ example, {\it Journal of Physics: Conference Series}, Vol 103,  DOI \url{https://iopscience.iop.org/article/10.1088/1742-6596/103/1/012006}, 
 and Quantum symmetries for exceptional $\SU{4}$ modular invariants associated with conformal embeddings,  {\it Symmetry, Integrability and Geometry : Methods and Applications}, { SIGMA 5} (2009), 044, 31 pp,  arXiv:0805.4678, \url{https://doi.org/10.3842/SIGMA.2009.044}
\bibitem{YellowBook} Di Francesco  P., Matthieu P. and  Senechal   D., Conformal field theory, Springer, (1997).
\bibitem{DiFrancescoZuber}  Di Francesco P. and   Zuber J.-B., $\SU{N}$ lattice integrable models associated with graphs,   {\it Nucl.  Phys.},  {B 338},  pp 602--646, (1990).  
\bibitem{Dorey:CoxeterElement} Dorey P., Partition Functions, Intertwiners and the Coxeter Element. arXiv:hep-th/9205040.  {\it Int. J. Mod. Phys} {A8}, pp 193-208, (1993).
\bibitem{EvansPughSU3} Evans D. E.  and Pugh M., Ocneanu cells and Boltzmann weights for the $\SU{3}$ ADE graphs. {\it M\"unster J. of Math.} {2},  pp 95-142 (2009)
\bibitem{Finkelberg} Finkelberg, M., An equivalence of fusion categories, {\it Geom. Funct. Anal.}  6 (1996), 249-267.
\bibitem{Huang} Y.-Z. Huang, Vertex operator algebras, the Verlinde conjecture, and modular tensor categories, {\it Proc. Natl. Acad. Sci. USA}, 102 (2005), 5352-5356.
\bibitem{Kac:book}  Kac V., Infinite dimensional Lie algebras, Cambridge University Press, Cambridge (1990).
\bibitem{KazhdanLusztig} Kazhdan D. and  Lusztig G., Tensor structures arising from affine Lie algebras, III,    {\it J. Amer. Math. Soc.},  {7}, pp 335--381, (1994).
\bibitem{KirilovOstrik} Kirillov A. and Ostrik V.,  On q-analog of McKay correspondence and ADE classification of SL2 conformal  field theories, {\it Adv. in Math.},  {171- 2}, pp 183--227, (2002).
\bibitem{Magma} Bosma W., Cannon J., and Playoust C., {\sl The Magma algebra system. I. The user language}, {\it J. Symbolic Comput.}, {\bf 24} (1997), 235--265, \url{http://magma.maths.usyd.edu.au}
\bibitem{Mathematica} Wolfram Research, Inc., Mathematica, Champaign, IL (2010).
\bibitem{Ocneanu:paths} Ocneanu  A., Paths on Coxeter diagrams: from Platonic solids and singularities to minimal models and subfactors,  Notes by Goto S.,  {\it Fields Institute Monographs}, Eds. Rajarama Bhat et al, (1999).  
\bibitem{Ocneanu:Bariloche}  Ocneanu  A., The Classification of  subgroups of quantum $\SU{N}$, in  ``Quantum symmetries in theoretical physics and mathematics'', Bariloche 2000, Eds. Coquereaux R., Garc\'{\i}a A. and Trinchero~R., {\it  AMS Contemporary Mathematics}, {294}, pp 133--160, (2000). 
\bibitem{Ocneanu:MSRI} Ocneanu  A., Higher Coxeter systems,  \url{http://www.msri.org/publications/ln/msri/2000/subfactors/ocneanu}, (2000). 
\bibitem{Ocneanu:posters} Ocneanu  A.,   Poster communications.
\bibitem{Ocneanu:Harvard_roots} A. Ocneanu. Harvard Lectures (2017-2018). YouTube: Video files  Adrian Ocneanu Harvard Physics L22,  267 2017 10 25, L23 , 267 2017 10 27,   L24, 267 2017 10 30, 
\url{https://www.youtube.com/watch?v=8ls_s7cpEjA&feature=youtu.be&t=2700}
\bibitem{OEIS} OEIS: The Online Encyclopedia of Integer Sequences, N.J.A. Sloane, \url{/https://oeis.org}
\bibitem{Ostrik}   Ostrik  V., Module categories, weak Hopf algebras and modular invariants, {\it Transform.  groups},  {8}, no 2, pp 177--206, (2003).
\bibitem{PleskenPohst} Plesken W. and Pohst M., Constructing integral lattices with prescribed minimum,  {\it Mathematics of Computation}, Vol 45, No 171, pp 209-221, and supplement S5-S16. 
\bibitem{Zagier:modularforms} Zagier D.B., Elliptic Modular Forms and Their Applications, in {\it  `The 1-2-3 of Modular forms'},  Lectures at a Summer School in Nordfjordeid, Norway, Springer (2008).
\bibitem{DezaGrishukhin} Deza M. and Grishukhin V., Delaunay Polytopes of Cut Lattices, {\it Linear Algebra and Its Applications}, 226-228:667-685 (1995).
\end{thebibliography}
\end{document}